\documentclass[reqno]{amsart}

\usepackage{amsmath,amssymb,amsthm,mathtools}
\usepackage{enumitem}
\usepackage{hyperref}
\usepackage{xcolor}
\usepackage[expansion=false]{microtype}
\usepackage{array}
\usepackage{booktabs}
\usepackage{tikz-cd}
\usepackage{todonotes}
\hypersetup{
colorlinks=true, linkcolor=blue!50!black, citecolor=blue!50!black, urlcolor=blue!50!black }

\newtheorem{theorem}{Theorem}[section]
\newtheorem{conjecture}[theorem]{Conjecture}
\newtheorem{proposition}[theorem]{Proposition}
\newtheorem{lemma}[theorem]{Lemma}
\newtheorem{corollary}[theorem]{Corollary}

\theoremstyle{definition}
\newtheorem{definition}[theorem]{Definition}
\theoremstyle{remark}
\newtheorem{remark}[theorem]{Remark}

\newcommand{\R}{\mathbb R}

\newcommand{\Z}{\mathbb Z}
\newcommand{\Pj}{\mathbb P}
\newcommand{\HH}{\mathbb H}

\newcommand{\Rad}{\operatorname{Rad}}
\newcommand{\rank}{\operatorname{rank}}
\newcommand{\Isom}{\mathrm{Isom}}
\newcommand{\SO}{\mathrm{SO}}
\newcommand{\SL}{\mathrm{SL}}
\newcommand{\SU}{\mathrm{SU}}
\newcommand{\Sp}{\mathrm{Sp}}

\newcommand{\PGL}{\mathrm{PGL}}
\newcommand{\PO}{\mathrm{PO}}
\newcommand{\Ogrp}{\mathrm{O}}
\newcommand{\Ein}{\mathrm{Ein}}
\newcommand{\Spin}{\mathrm{Spin}}

\newcommand{\Ad}{\operatorname{Ad}}
\newcommand{\id}{\operatorname{id}}

\newcommand{\Zcl}{\operatorname{Zcl}}
\newcommand{\limcone}{\mathcal L}

\newcommand{\ip}[2]{\left\langle #1,#2\right\rangle}

\title[$C^1$-limit-set rigidity]{$C^1$-rigidity of $\HH^{p,q}$-convex cocompact representations of Gromov hyperbolic groups}
\author{Nikolay Bogachev}
\address{Department of Computer and Mathematical Sciences, University of Toronto Scarborough, 1095 Military Trail, Toronto, ON M1C 1A3, Canada}
\email{n.bogachev@utoronto.ca}

\begin{document}

\begin{abstract}
Let $p>3$ and let $\Gamma$ be a non-elementary Gromov hyperbolic group with $\partial_\infty\Gamma\cong\mathbb S^{p-1}$. We prove that every $\HH^{p,2}$-convex cocompact representation
$$
\rho:\Gamma\longrightarrow \PO_{p,3}
$$
whose proximal limit set $\Lambda_\rho\subset\partial_\infty\HH^{p,2}$ is a $C^1$ submanifold preserves a totally geodesic copy of $\HH^p\subset\HH^{p,2}$. Consequently, $\Gamma/\ker \rho$ is virtually a cocompact lattice in $\Isom(\HH^p)$.

For odd $p$, we prove the same conclusion for every signature $(p,q)$ with $p>q\ge2$, under the same hypothesis $\partial_\infty\Gamma\cong\mathbb S^{p-1}$.

Our methods combine dynamics, differential topology, and representation theory. The key new ingredient is to study the radicals of the semidefinite forms induced on the tangent spaces of the limit set. They form a continuous vector bundle of constant rank. The tangent action at attracting fixed points gives a bound, in terms of the radical rank, on the number of eigenvalues of an orthogonal lift whose modulus exceeds one. We then use the topology of the sphere and the representation theory of real algebraic groups to exclude the remaining possibilities. In arbitrary signature with $p>q$, these arguments already imply that a representation with a $C^1$ limit set cannot be Zariski dense.
\end{abstract}

\maketitle

\centerline{{\itshape Dedicated to the memory of Michael Kapovich (1963--2026)}}

\section{Introduction}

Classical quasi-Fuchsian groups are flexible, but regularity of their limit sets is rigid: sufficiently regular limit sets typically force the group to be Fuchsian. This principle has many avatars. In anti-de Sitter geometry it was studied by Barbot--M\'erigot \cite{BM12}, and Glorieux--Monclair \cite[Theorem~1.3]{GM24} proved that an AdS quasi-Fuchsian group whose limit set is $C^1$ is Fuchsian. In higher-rank Anosov geometry, regularity of the limit set imposes strong spectral restrictions; see Zimmer \cite[Theorem~1.46]{Zim21}, Zhang--Zimmer \cite[Theorems~1.12 and~1.14]{ZZ24}, and, at the $C^1$ threshold, Canary--Oh--Zimmer \cite[Theorem~1.5]{COZ24}. The Lipschitz regime was studied systematically by Pozzetti--Sambarino--Wienhard \cite{PSW23}; their Proposition~8.2 includes the graph description relevant to pseudo-hyperbolic convex cocompactness.

We write $\HH^{p,q-1}$ for the pseudo-hyperbolic space whose projective isometry group is $\PO_{p,q}$. Limit sets of $\HH^{p,q-1}$-convex cocompact representations of groups with $\partial_\infty\Gamma\cong\mathbb S^{p-1}$ are always Lipschitz submanifolds, see \cite[Proposition~8.2]{PSW23}. We prove that, in the signatures considered below, they can be $C^1$ only when the representation preserves a totally geodesic real hyperbolic space. Thus regularity becomes rigid already at $C^1$, although the Lipschitz condition holds throughout the deformation space. 

The problem belongs to higher higher Teichm\"uller theory. Beyrer--Kassel \cite[Theorem~1.3]{BK25} prove that, for a Gromov hyperbolic group of virtual cohomological dimension $p$, pseudo-hyperbolic convex cocompact representations into $\PO_{p,q}$ form a union of connected components. In particular, the component containing the standard representation of a closed hyperbolic $p$-manifold group consists entirely of such representations \cite[Corollary~1.4]{BK25}. Our results identify the representations with $C^1$ limit sets in these components whenever $q=3$ or $p$ is odd, with $p>q$.

We address the following conjecture.

\begin{conjecture}[General $C^1$-rigidity]\label{conj:pq}

Let $p>q\ge2$, let $\Gamma$ be a non-elementary Gromov hyperbolic group with $\partial_\infty\Gamma\cong\mathbb S^{p-1}$, and let $\rho:\Gamma\to\PO_{p,q}$ be $\HH^{p,q-1}$-convex cocompact. If the proximal limit set $\Lambda_\rho\subset\partial_\infty\HH^{p,q-1}$ is a $C^1$ submanifold, then $\rho(\Gamma)$ preserves a totally geodesic copy of $\HH^p$ and therefore $\Gamma/\ker \rho$ is virtually a uniform hyperbolic lattice.

\end{conjecture}

For $q=2$, this is the theorem of Glorieux--Monclair \cite[Theorem~1.3]{GM24}, following the work of Barbot--M\'erigot \cite{BM12} on AdS quasi-Fuchsian representations. We prove the conjecture for $q=3$ and arbitrary $p>3$, and for every $p>q\ge2$ with $p$ odd. The case $q=3$ is the first beyond anti-de Sitter geometry.

\begin{theorem}[Index-three $C^1$-rigidity]\label{thm:main}
Let $p>3$, let $\Gamma$ be a non-elementary Gromov hyperbolic group with $\partial_\infty\Gamma\cong\mathbb S^{p-1}$, and let $\rho:\Gamma\to\PO_{p,3}$ be an $\HH^{p,2}$-convex cocompact representation. If the proximal limit set $\Lambda_\rho\subset\partial_\infty\HH^{p,2}$ is a $C^1$ submanifold, then $\rho(\Gamma)$ preserves a totally geodesic copy of $\HH^p\subset\HH^{p,2}$. Consequently, $\Gamma/\ker\rho$ is virtually a cocompact lattice in $\Isom(\HH^p)$.

\end{theorem}

Our second result treats all negative indices under a parity assumption on $p$.

\begin{theorem}[Odd-$p$ $C^1$-rigidity in general signature]\label{thm:oddpq}

Let $p>q\ge2$, with $p$ odd, and let $\Gamma$ be a non-elementary Gromov hyperbolic group such that $\partial_\infty\Gamma\cong\mathbb S^{p-1}.$ Let $\rho:\Gamma\to\PO_{p,q}$ be an $\HH^{p,q-1}$-convex cocompact representation. If its proximal limit set $\Lambda_\rho$ is a $C^1$ submanifold, then $\rho(\Gamma)$ preserves a totally geodesic copy of $\HH^p \subset \HH^{p,q-1}$. Equivalently, after conjugacy, $\rho(\Gamma)$ is contained in the image of $\Ogrp_{p,1}\times\Ogrp_{q-1}$ in $\PO_{p,q}$. Moreover, the projection of $\rho(\Gamma)$ to $\PO_{p,1}$ is a uniform lattice; in particular, $\Gamma/\ker\rho$ is virtually a uniform lattice in $\operatorname{Isom}(\HH^p)$.

\end{theorem}

The proof of Theorem~\ref{thm:main} is considerably more involved than that of Theorem~\ref{thm:oddpq}. For odd $p$, an Euler-class argument makes the limit set spacelike, after which entropy rigidity applies. The index-three theorem requires a detailed analysis of the possible Zariski closures; we explain this difference below.

We now compare these results with regularity rigidity results in other settings. The $C^1$ threshold contrasts with the regularity of the boundaries of strictly convex divisible domains. The corresponding holonomy representations form Benoist components \cite{Ben05}. The boundaries of these domains are always $C^1$, while $C^2$ regularity forces the domain to be an ellipsoid; see Benoist \cite[Theorems~1.1 and~1.3]{Ben04}. Between these thresholds, the optimal H\"older regularity of the tangent hyperplane map is determined by eigenvalue ratios, or equivalently by the limit cone, as shown by Guichard \cite{Gui05}. This applies in particular to the convex projective domains associated with Hitchin representations in $\SL_3(\R)$. Related regularity rigidity results include Potrie--Sambarino's theorem \cite[Theorem~D]{PS17} that a Hitchin representation with a smooth Frenet curve is Fuchsian, and the results of Zimmer \cite[Theorem~1.46]{Zim21} and Zhang--Zimmer \cite[Theorems~1.12 and~1.14]{ZZ24}.

Here every limit set is a Lipschitz $(p-1)$-sphere, so its Hausdorff dimension, measured using an auxiliary Riemannian metric, is always $p-1$. This also differs from the behavior near the reducible block embedding of a Fuchsian surface group into $\SL_3(\R)$: Li--Pan--Xu \cite[Theorem~1.3]{LPX23} prove that, as irreducible representations converge to this embedding, the Hausdorff dimensions of their limit sets tend to $3/2$, whereas the original projective line has dimension $1$. In our setting, loss of $C^1$ regularity is invisible to Hausdorff dimension.

Regularity of the limit set as an embedded submanifold should also be distinguished from regularity of an equivariant map between two limit sets. The optimal H\"older exponents of such maps can vary even when both limit sets are smooth. For Fuchsian representations, Thurston's asymmetric metric can be expressed in terms of these exponents, and related constructions give asymmetric metrics on Hitchin components; see Nolte \cite[Theorem~A and Sections~1.1--1.2]{Nol24b}. In another direction, Canary--Oh--Zimmer \cite[Theorem~1.1]{COZ24} prove that, for Zariski-dense discrete subgroups of simple real algebraic groups and maximal parabolic quotients, a bi-Lipschitz equivariant map between limit sets extends to a diffeomorphism of the ambient flag varieties. Our theorems concern the embedded limit set itself and impose no regularity assumption on a parametrized boundary map.

We now turn to the main ideas of the proofs. Our approach is to study the radicals of the forms induced on the tangent spaces of the limit set. In the AdS case, Glorieux--Monclair \cite{GM24} first find a spacelike tangent point by viewing the limit set as the graph of a scalar-valued function and taking a maximum. Minimality then makes every tangent space spacelike; the derivative at a fixed point gives a spectral restriction, and their description of the possible Zariski closures completes the proof; see \cite[Proposition~3.3, Lemma~4.1, Corollary~4.2, and Lemma~4.3]{GM24}. The maximum argument does not extend to a map with higher-dimensional target, as they explain after Corollary~4.2. We study the possibly nonzero tangent radicals and show that they have constant rank and form a continuous vector bundle. When $q=3$, the intermediate radical rank $r=1$ must also be excluded, and this accounts for much of the proof.

Passing from a geometric condition to restrictions on eigenvalues and then using the limit cone is a well-known approach to rigidity. Examples include Benoist's work on boundary regularity \cite[Section~6]{Ben04}, Dal'Bo--Kim's marked-length rigidity theorem \cite[Section~2]{DK02}, the regularity estimates of Zhang--Zimmer \cite{ZZ24}, and Nolte's study of leaves of foliated projective structures \cite[Section~4.3]{Nol24a}. Our contribution is the use of the tangent radical to obtain the spectral restriction when spacelikeness is not known, and the arguments needed to exclude the resulting algebraic possibilities.

The main algebraic result is the following.

\begin{theorem}\label{thm:intro-smallwall}

Let $H<\SO_{p,3}$ be connected, noncompact, reductive, and irreducible on $\R^{p,3}$, with $p>3$. Then the ambient limit cone of $H$ is not contained in the wall $\{\lambda_3=0\}$.

\end{theorem}

Theorem~\ref{thm:intro-smallwall} is related to the classification of connected irreducible subgroups of $\SO_{2,n}$ by Di Scala--Leistner \cite{DSL11}, which Glorieux--Monclair use in the AdS case \cite[Proposition~3.3]{GM24}. We only classify the representations compatible with the indicated wall. The wall condition bounds the total dimension of the nonzero restricted-weight spaces, but gives no direct bound on the zero-weight space. One must therefore determine the possible irreducible real modules and their invariant signatures, taking account of real, complex, and quaternionic representation types. General results constrain Zariski closures in the divisible convex setting \cite[Theorem~1.1]{Ben03} and for positive representations in split groups \cite[Theorem~A]{Sam24}. Their hypotheses do not apply directly to the negative projective limit spheres considered here, so the required algebraic analysis is carried out in Section~\ref{sec:smallwall}.

\subsection*{Line of argument and organization of the paper}

Section~\ref{sec:prelim} recalls the geometry and fixes the conventions for Cartan and Jordan projections. Each component of the lift of the limit set to the double cover $\mathbb S^{p-1}\times\mathbb S^{q-1}$ is the graph of a strictly $1$-Lipschitz map $\varphi:\mathbb S^{p-1}\longrightarrow\mathbb S^{q-1}.$ If the limit set is $C^1$, then so is $\varphi$ by Lemma~\ref{lem:c1graph}. The induced form
$$
Q=g_{\mathbb S^{p-1}}-\varphi^*g_{\mathbb S^{q-1}}
$$
is continuous and positive semidefinite. Minimality of the boundary action forces its radical to have constant rank $r\in\{0,\ldots,q-1\}$; continuity of $Q$ then makes the radical a continuous subbundle, as proved in Proposition~\ref{prop:constant}.

Section~\ref{sec:r0r2} treats the two extreme ranks. If $r=0$, the limit set is spacelike, and rigidity follows from Monclair's theorem \cite[Theorem~3.11]{Mon25}. We recall a proof using the critical exponent and entropy rigidity. If $r=q-1$, the map $\varphi$ would be a submersion and hence surjective. This is impossible for a strictly distance-decreasing map between unit spheres.

Section~\ref{sec:wall} examines the derivative action at attracting fixed points. A radical of rank $r\le q-2$ forces the limit cone into the face $F_{r+1}=\{\lambda_{r+2}=0\}$ of the ambient closed Weyl chamber. In particular, the image cannot be Zariski dense. The same section excludes invariant totally isotropic subspaces and proves reductivity of the identity component of the Zariski closure. Section~\ref{sec:descent} shows that, when the Zariski closure acts reducibly, the limit set spans a nondegenerate subspace of signature $(p,b)$ with $b<q$. Restriction to this subspace reduces the problem to a smaller negative index.

For $q=3$, both $r=0$ and $r=1$ give the wall condition $\lambda_3=0$. Section~\ref{sec:smallwall} proves Theorem~\ref{thm:intro-smallwall} by reducing to a noncompact simple group of real rank at most two and classifying its orthogonal modules with the required weight bound. Section~\ref{sec:proofmain} then completes the proof of Theorem~\ref{thm:main}; the reducible case reduces to the theorem of Glorieux--Monclair \cite[Theorem~1.3]{GM24}. For odd $p$, Section~\ref{sec:pqodd} uses the nonzero Euler class of $T\mathbb S^{p-1}$ to rule out every radical subbundle of positive rank, proving Theorem~\ref{thm:oddpq}. The last section discusses the remaining even-$p$ cases with $q>3$.

\subsection*{Funding}
This research was supported by NSERC Discovery Grant RGPIN-2024-05680.

\subsection*{Acknowledgements}
I first learned about this problem from Olivier Glorieux during my visit to IH\'ES in December 2019 -- February 2020. I am grateful to him for introducing me to the problem, for explaining his joint work \cite{GM24} with Daniel Monclair, and for pointing out the main difficulties. I also thank Alex Nolte, whom I met at IH\'ES in the summer of 2025, for fruitful discussions, helpful comments on this problem and on my approach to it, as well as his remarks and suggestions on the early version of this manuscript. I am grateful to Fanny Kassel and to IH\'ES for their hospitality and for the stimulating atmosphere during my visits. I am especially grateful to the late Misha Kapovich for our many conversations, mathematical and otherwise, and for all that I learned from him about applying ideas from dynamics and topology to problems in group theory.

\subsection*{Use of AI}
The starting point for this paper was Beyrer--Kassel's description of negative limit spheres as graphs of strictly $1$-Lipschitz maps \cite[Lemma~3.3 and Corollary~3.4]{BK25}. It was then natural to study the tangent radical whose rank measures how far a $C^1$ limit set can be from being spacelike. I used large language models while developing the paper: to discuss mathematical ideas, explore proof strategies, check calculations and proofs, find references, and write Python code that I initially used for the classification of low-rank algebraic groups. The proof presented here does not depend on this code. I also used these models to help formulate and revise parts of my exposition. I take responsibility for the mathematical claims and proofs.

\section{Pseudo-hyperbolic geometry and the tangent radical}\label{sec:prelim}

\subsection{The projective model and the Einstein boundary}

Throughout the paper, we mostly follow the notation and conventions of \cite{BK25}; see also \cite{DGK18} for the original definitions. All round spheres are unit spheres; thus $\mathbb S^m$ carries its standard metric of sectional curvature $1$. If $W$ is a Euclidean vector space, $\mathbb S(W)$ denotes its unit sphere.

Let $a\ge2$ and $b\ge0$, and let $V=\R^{a,b+1}$ be endowed with a nondegenerate symmetric bilinear form $b_V$ of signature $(a,b+1)$. The pseudo-hyperbolic space of signature $(a,b)$ is
$$
\HH^{a,b}=\{[v]\in\Pj(V) \mid b_V(v,v)<0\}.
$$

It is double-covered by the unit pseudo-hyperboloid
$$
\widehat{\HH}^{a,b}=\{x\in V \mid b_V(x,x)=-1\},
$$
with deck transformation $x\mapsto -x$. Note that $\widehat{\HH}^{a,b}$ is homeomorphic to $\R^a\times\mathbb S^b$ and therefore is connected for $b>0$, but has two connected components when $b=0$ ($\widehat{\HH}^{a,0}$ is a two-sheeted hyperboloid). At $x\in\widehat{\HH}^{a,b}$ the tangent space is
$$
T_x\widehat{\HH}^{a,b}=x^\perp =\{y\in V\mid b_V(x,y)=0\},
$$
where the restriction of $b_V$ to $x^\perp$ has signature $(a,b)$. It defines a pseudo-Riemannian metric of constant sectional curvature $-1$. Its full projective isometry group is $\PO_{a,b+1}$. Its unparameterized geodesics are the connected components of its intersections with projective lines. Equivalently, $\HH^{a,b}$ is the pseudo-Riemannian symmetric space $\PO_{a,b+1}/\Ogrp_{a,b}$. For $b=0$ this is real hyperbolic space $\HH^a$, while $\HH^{a,1}$ is anti-de Sitter space. See \cite[Section~2.2]{BK25} and \cite[Sections~1.1 and~2.1]{DGK18}.

If $W<V$ is nondegenerate of signature $(c,d+1)$, then
$$
\HH(W)=\Pj(W)\cap\HH^{a,b}
$$
is a totally geodesic copy of $\HH^{c,d}$. In particular, the copy of $\HH^p$ in Theorem~\ref{thm:main} is determined by a subspace of signature $(p,1)$.

The ideal boundary in the projective model is the quadric
$$
\partial_\infty\HH^{a,b} = \{[v]\in\Pj(V)\mid b_V(v,v)=0\}.
$$

For an isotropic line $\ell<V$, there is a canonical identification
$$
T_\ell\partial_\infty\HH^{a,b} \cong\operatorname{Hom}(\ell,\ell^\perp/\ell).
$$

The form induced by $b_V$ on $\ell^\perp/\ell$ has signature $(a-1,b)$. Choosing a nonzero vector in $\ell$ therefore gives a quadratic form of signature $(a-1,b)$ on $T_\ell\partial_\infty\HH^{a,b}$; rescaling that vector changes the form by a positive factor. In this way the boundary inherits a conformal pseudo-Riemannian structure and is in fact the {\em Einstein universe} $\partial_\infty\HH^{a,b}=\Ein^{a-1,b}.$

For the projectivized-null-cone definition and its conformal structure, see \cite[Section~2.1]{FM10}; see also \cite[Sections~2.1--2.2]{BK25} for the boundary of $\HH^{a,b}$ in these conventions.

The double cover of the Einstein universe admits a useful product model. Namely, the space of positive isotropic rays is
$$
\partial_\infty\widehat{\HH}^{a,b} =\bigl(\{v\ne0 \mid b_V(v,v)=0\}\bigr)/\R_{>0}.
$$

Choose an orthogonal splitting $V=V_+\oplus V_-$ with $b_V$ positive on $V_+$, negative on $V_-$, and $\dim V_+=a$, $\dim V_-=b+1$. Every positive isotropic ray has a unique representative $(u,w)$ with respect to the Euclidean forms $b_V\mid_{V_+}$ and $-b_V\mid_{V_-}$, with
$$
u\in\mathbb S(V_+)\cong\mathbb S^{a-1},\qquad w\in\mathbb S(V_-)\cong\mathbb S^b.
$$

Consequently, there are diffeomorphisms
$$
\partial_\infty\widehat{\HH}^{a,b} \cong\mathbb S^{a-1}\times\mathbb S^b, \qquad \partial_\infty\HH^{a,b} = \Ein^{a-1,b} \cong(\mathbb S^{a-1}\times\mathbb S^b)/\{\pm1\},
$$
where the deck transformation is $(u,w)\mapsto(-u,-w)$. On the double cover the Einstein conformal class is represented by $g_{\mathbb S^{a-1}}-g_{\mathbb S^b}.$ See \cite[Section~2.1]{FM10} for the conformal double cover and \cite[Proposition~2.3, Lemmas~2.5--2.6, and~(2.7)]{BK25} for this explicit sphere-product model. The identification with a product depends on the splitting $V_+\oplus V_-$, whereas the conformal structure and the projective quotient are intrinsic.

\subsection{Cartan and Jordan projections and the limit cone}\label{sec:projections}

We recall the conventions that will be used below. Let $G$ be a real reductive group, let $\mathfrak g=\mathfrak k\oplus\mathfrak p$ be a Cartan decomposition of its Lie algebra, and let $K<G$ be the corresponding maximal compact subgroup. Choose a maximal abelian subspace $\mathfrak a_G\subset\mathfrak p$. We use the convention that $\mathfrak a_G^+$ denotes a {\em closed positive Weyl chamber} in $\mathfrak a_G$. The (possibly extended) Weyl group $W_G=N_G(\mathfrak a_G)/Z_G(\mathfrak a_G)$ acts on $\mathfrak a_G$, with $\mathfrak a_G^+$ as a fundamental domain; here $N_G$ and $Z_G$ denote the normalizer and centralizer in $G$ via the adjoint action, respectively. The Cartan decomposition
$$
G=K\exp(\mathfrak a_G^+)K
$$
associates with every $g\in G$ a unique element $\mu_G(g)\in\mathfrak a_G^+$, called its \emph{Cartan projection}. Thus $g=k\exp(\mu_G(g))k'$ for some $k,k'\in K$.

There is also a map $\lambda_G:G\longrightarrow\mathfrak a_G^+$, called the \emph{Jordan projection}, or \emph{Lyapunov projection}, which can be defined as follows. By the real Jordan decomposition, every $g\in G$ can be written uniquely as a commuting product $g=g_hg_eg_u$, where $g_h, g_e, g_u\in G$ are, respectively, hyperbolic, elliptic, and unipotent. Here $g_h$ is conjugate to an element of $\exp(\mathfrak a_G^+)$, while $g_e$ is conjugate to an element of $K$. We denote by $\lambda_G(g)$ the unique element of $\mathfrak a_G^+$ such that $g_h$ is conjugate to $\exp(\lambda_G(g))$. The Cartan and Jordan projections are related by
$$
\lambda_G(g)=\lim_{n\longrightarrow\infty}\frac1n\mu_G(g^n).
$$

Informally, the Cartan projection records singular values, whereas the Jordan projection records eigenvalue moduli and their asymptotic growth. We refer to \cite[Sections~2.3--2.4]{GGKW17} for these definitions and to \cite[Section~2.5]{BK25} for the Jordan projection in the conventions used here.

For a subsemigroup $\Delta\subset G$, its \emph{Benoist limit cone} is
$$
\limcone^G_\Delta
=\overline{\{t\lambda_G(\delta):t\ge0,\ \delta\in\Delta\}}
\subset\mathfrak a_G^+.
$$

The superscript will usually be omitted when the ambient group is clear. If $G$ is connected and semisimple and $\Delta$ is Zariski dense in $G$, this cone is convex and has nonempty interior in $\mathfrak a_G$ by Benoist's theorem \cite[Theorem~1.2]{Ben97}.

These notions depend on the ambient group. In particular, if $L<G$ is reductive, let $A_L<L$ be a maximal real split torus and put $\mathfrak a_L=\operatorname{Lie}(A_L)$. Its intrinsic Cartan and Jordan projections take values in a closed chamber $\mathfrak a_L^+$ for $L$, while the ambient projections take values in $\mathfrak a_G^+$, and the two should not be identified. After conjugating $A_L$ into $A_G:=\exp(\mathfrak a_G)$, we regard $\mathfrak a_L$ as a subspace of $\mathfrak a_G$. The ambient Jordan projection of $\exp X$, for $X\in\mathfrak a_L$, is then the unique $W_G$-translate of $X$ lying in $\mathfrak a_G^+$. All Weyl chambers, coordinates, walls, and faces appearing below are those of the ambient orthogonal group unless explicitly stated otherwise.

In particular, the set of ambient Jordan projections of elements of $L$ is
$$
\lambda_G(L)=\bigcup_{w\in W_G}\bigl(w\mathfrak a_L\cap\mathfrak a_G^+\bigr),
$$
and this union need not be the intersection of a single linear subspace with $\mathfrak a_G^+$. For example, take $L=\operatorname{diag}(\SL_3(\R),1)<\SL_4(\R)$. If $a\ge b\ge c$ and $a+b+c=0$, the ambient Jordan projection sorts the four coordinates $(a,b,c,0)$, giving $(a,b,0,c)$ when $b\ge0$ and $(a,0,b,c)$ when $b\le0$. The image is the union of two planar sectors. It contains $(2,1,0,-3)$ and $(3,0,-1,-2)$, but not their midpoint, which has no zero coordinate. It is therefore not convex. This example is discussed for Cartan projections by Dey--Oh \cite[Section~7]{DO25}; the same calculation gives the Jordan image.

We now describe these objects for $G=\PO_{a,b+1}$. Let $V=\mathbb R^{a,b+1}$, endowed with its defining bilinear form $b_V$. For the standard orthogonal splitting of $V$ into positive- and negative-definite summands, take $K=\mathrm P(\Ogrp_a\times\Ogrp_{b+1})$. When referring to singular values, we use the corresponding Euclidean scalar product, obtained by reversing the sign of $b_V$ on the negative-definite summand. Put $s=\min(a,b+1)$ and $d=|a-(b+1)|$. Choose a Witt basis $u_1,\ldots,u_s,z_1,\ldots,z_d,v_1,\ldots,v_s$ of $V$, where the spans of the $u_i$ and of the $v_i$ are totally isotropic, $b_V(u_i,v_j)=\delta_{ij}$, and the span of the $z_j$ is definite and orthogonal to all the $u_i$ and $v_i$. We choose this basis to be orthonormal for the preceding Euclidean scalar product; in other words, it is a $K$-adapted Witt basis.

For $\mathbf t=(t_1,\ldots,t_s)\in\mathbb R^s$, set
$$
H(\mathbf t) =\operatorname{diag}\left(t_1,\ldots,t_s, \underbrace{0,\ldots,0}_{d},-t_1,\ldots,-t_s\right),
\qquad
\mathfrak a_G=\{H(\mathbf t):\mathbf t\in\mathbb R^s\}.
$$

Thus, $\mathfrak a_G$ is the Lie algebra of a maximal real split torus of $G$. We identify $\mathfrak a_G$ with $\mathbb R^s$ by sending $H(\mathbf t)$ to $\mathbf t$. If $a\ne b+1$, the restricted root system is of type $B_s$, and under this identification we take
$$
\mathfrak a_G^+ =\{\mathbf t\in\mathbb R^s:t_1\ge\cdots\ge t_s\ge0\}.
$$

When $a=b+1=s$, the restricted root system of the identity component $G^\circ$ is of type $D_s$. The usual closed chamber for $G^\circ$ and the folded chamber for the full group $G=\PO_{s,s}$ are, respectively,
$$
t_1\ge\cdots\ge t_{s-1}\ge|t_s| \qquad\text{and}\qquad t_1\ge\cdots\ge t_s\ge0.
$$

In the first chamber $t_s$ is signed, whereas the last coordinate in the second is $|t_s|$. This distinction will not enter our arguments, since all the wall statements below concern $\PO_{p,q}$ with $p>q$.

For $g\in G$, write $\mu_G(g)=(\mu_1(g),\ldots,\mu_s(g))$ and $\lambda_G(g)=(\lambda_1(g),\ldots,\lambda_s(g))$. If $\widetilde g\in\Ogrp_{a,b+1}$ is any lift of $g$, then the logarithms of its singular values, counted with multiplicity, and the logarithms of the moduli of its eigenvalues, counted with algebraic multiplicity, are, respectively,
$$
\begin{aligned}
&\mu_1(g),\ldots,\mu_s(g),
\underbrace{0,\ldots,0}_{d},
-\mu_s(g),\ldots,-\mu_1(g),\\
&\lambda_1(g),\ldots,\lambda_s(g),
\underbrace{0,\ldots,0}_{d},
-\lambda_s(g),\ldots,-\lambda_1(g).
\end{aligned}
$$

Here $\mu_1(g)\ge\cdots\ge\mu_s(g)\ge0$ and $\lambda_1(g)\ge\cdots\ge\lambda_s(g)\ge0$. Replacing the lift $\widetilde g$ by $-\widetilde g$ changes neither projection. See \cite[Example~2.9]{BK25} for the Jordan-projection description in the unbalanced case, which is the case used below.

For $G=\PO_{p,q}$ with $p>q$, the real rank is $q$, and
$$
\mathfrak a_G^+ =
\{(t_1,\ldots,t_q)\in\mathbb R^q:t_1\ge\cdots\ge t_q\ge0\}.
$$

The restricted root system is of type $B_q$, with simple roots $\alpha_i(\mathbf t)=t_i-t_{i+1}$ for $1\le i<q$ and $\alpha_q(\mathbf t)=t_q$. The simple-root wall corresponding to $\alpha_i$ is $\mathfrak a_G^+\cap\ker\alpha_i$. For $1\le m\le q-1$, define
$$
\begin{aligned}
F_m
&:=\mathfrak a_G^+\cap\{\mathbf t\in\mathbb R^q:t_{m+1}=0\}\\
&=\mathfrak a_G^+\cap\bigcap_{j=m+1}^q\ker\alpha_j\\
&=\{(t_1,\ldots,t_m,0,\ldots,0):t_1\ge\cdots\ge t_m\ge0\}.
\end{aligned}
$$

Thus, $\dim F_m=m$, and $F_m$ is the codimension-$(q-m)$ face obtained by intersecting the last $q-m$ simple-root walls. In particular, $F_{q-1}=\mathfrak a_G^+\cap\{t_q=0\}$ is the last Weyl wall. The coordinates of $\mu_G$ are the logarithms of the largest $q$ singular values, while those of $\lambda_G$ are the logarithms of the largest $q$ eigenvalue moduli. The simple root corresponding to the stabilizer of an isotropic line is $\alpha_1$; the associated Cartan and Jordan gaps are $\mu_1-\mu_2$ and $\lambda_1-\lambda_2$, respectively.

\subsection{Convex cocompactness and \texorpdfstring{$P_1$}{P1}-Anosov representations}

Convexity in $\HH^{a,b}$ is projective convexity. More precisely, a subset of $\Pj(V)$ is {\em convex} if it is contained in an affine chart and is convex there; it is {\em properly convex} if its closure is contained in an affine chart and is convex there. These are the conventions of \cite[Section~1.1]{DGK18}. This notion should not be confused with geodesic convexity in the associated Riemannian symmetric space $X_{a,b+1}=\PO_{a,b+1}/K$, whose points are the negative-definite $(b+1)$-planes in $V$. By contrast, $\HH^{a,b}$ consists of negative lines and carries an indefinite metric when $b>0$. The two spaces coincide only for $b=0$.

\begin{definition}\label{def:hcc}

An infinite discrete subgroup $\Delta<\PO_{a,b+1}$ is \emph{$\HH^{a,b}$-convex cocompact} if it acts properly discontinuously and cocompactly on a closed properly convex subset $\mathcal C\subset\HH^{a,b}$ with nonempty interior, and if its ideal boundary
$$
\partial_i\mathcal C
=\overline{\mathcal C}^{\,\Pj(V)}\setminus\mathcal C =\overline{\mathcal C}^{\,\Pj(V)}\cap\partial_\infty\HH^{a,b}
$$
contains no nontrivial projective segment. A representation into $\PO_{a,b+1}$ is said to be $\HH^{a,b}$-convex cocompact if it has finite kernel and an image that is discrete and $\HH^{a,b}$-convex cocompact.
\end{definition}

This is \cite[Definition~1.2]{DGK18}; for $b\ge1$, see also \cite[Definition~1.2]{BK25}.

For $b=0$, this recovers the classical rank-one notion. For the Anosov characterization below, assume that $b\ge1$. Let $P_1<\PO_{a,b+1}$ be the stabilizer of an isotropic line. Then $\PO_{a,b+1}/P_1=\partial_\infty\HH^{a,b}$. We use the folded spectral coordinates introduced in Subsection~\ref{sec:projections}; thus $\lambda_1(g)$ and $\lambda_2(g)$ are the two largest nonnegative logarithms of eigenvalue moduli, and the $P_1$ gap is $\lambda_1(g)-\lambda_2(g)$.

\begin{definition}\label{def:p1anosov}
Let $\Gamma$ be a Gromov hyperbolic group. A representation $\rho:\Gamma\to\PO_{a,b+1}$ is \emph{$P_1$-Anosov} if there is a continuous $\rho$-equivariant map
$$
\xi:\partial_\infty\Gamma\longrightarrow\partial_\infty\HH^{a,b}
$$
with the following properties:

\begin{enumerate}[label=\textup{(\roman*)}]
\item $\xi$ is transverse, namely $\xi(\eta)\not\subset\xi(\eta')^\perp$ whenever $\eta\ne\eta'$; \item $\xi$ is dynamics-preserving: if $\gamma\in\Gamma$ has infinite order, then $\xi(\gamma^+)$ is the attracting isotropic line of $\rho(\gamma)$, where $\gamma^+$ is the attracting fixed point of $\gamma$ in $\partial_\infty\Gamma$; \item for every sequence $(\gamma_n)$ in pairwise distinct conjugacy classes,
$$
(\lambda_1-\lambda_2)(\rho(\gamma_n))\longrightarrow+\infty.
$$
\end{enumerate}
\end{definition}

This eigenvalue-gap characterization is \cite[Theorem~4.2]{GGKW17}; in the present orthogonal setting it is stated as \cite[Definition~2.17]{BK25}. An element $g\in\PGL(V)$ is \emph{proximal} if a lift has a simple real eigenvalue of strictly maximal modulus. Write $\ell_g^+$ for its attracting eigenline. We denote by $\Lambda_\rho$ the \emph{proximal limit set} of $\rho(\Gamma)$, defined by $\Lambda_\rho:=\overline{\{\ell_g^+:g\in\rho(\Gamma)\text{ is proximal}\}}\subset\Pj(V)$; see \cite[Definition~2.3 and Remark~2.4]{DGK18}.

If $\rho$ is $P_1$-Anosov, then $\rho(\gamma)$ is proximal for every infinite-order $\gamma\in\Gamma$, and the dynamics-preserving property gives $\ell_{\rho(\gamma)}^+=\xi(\gamma^+)$. Indeed, applying condition~\textup{(iii)} to the positive powers of $\gamma$ gives $(\lambda_1-\lambda_2)(\rho(\gamma))>0$. Since the attracting fixed points $\gamma^+$ are dense in $\partial_\infty\Gamma$, it follows that $\Lambda_\rho=\xi(\partial_\infty\Gamma)$.

A subset of $\partial_\infty\HH^{a,b}$ is \emph{negative} if it has a lift to the isotropic cone on which $b_V(v,w)<0$ for every pair of nonproportional lifted vectors. The terminology is \cite[Definition~2.7]{BK25}. In the Lorentzian case $b=1$, negative sets are traditionally called acausal; in higher index we use ``negative.''

\begin{theorem}[Danciger--Gu\'eritaud--Kassel]\label{thm:dgk-characterization} 
Let $a\ge2$ and $b\ge1$, let $\Gamma$ be an infinite Gromov hyperbolic group, and let $\rho:\Gamma\longrightarrow\PO_{a,b+1}$ be a representation. The following conditions are equivalent:

\begin{enumerate}[label=\textup{(\roman*)}]
\item $\rho$ is $\HH^{a,b}$-convex cocompact; 
\item $\rho$ is $P_1$-Anosov and $\rho(\Gamma)$ preserves a nonempty properly convex open subset of $\HH^{a,b}$; 
\item $\rho$ is $P_1$-Anosov and its proximal limit set $\Lambda_\rho\subset\partial_\infty\HH^{a,b}$ is negative.
\end{enumerate}

When these conditions hold, every invariant convex set $\mathcal C$ appearing in Definition~\ref{def:hcc} has ideal boundary $\partial_i\mathcal C=\Lambda_\rho=\xi(\partial_\infty\Gamma)$.
\end{theorem}

Applied to the discrete image $\Delta=\rho(\Gamma)$, this follows from \cite[Fact~2.18]{BK25}. The passage from the subgroup formulation to the representation formulation uses the fact that a $P_1$-Anosov representation has finite kernel and discrete image; see \cite[Theorem~1.3 and Remark~1.5]{GGKW17}. For the original subgroup formulation, see in particular \cite[Theorem~1.24]{DGK24}; the equivalence between $\HH^{a,b}$-convex cocompactness and the negative $P_1$-Anosov condition is the equivalence between conditions \textup{$(3^-)$} and \textup{$(5^-)$} there. Recall also that a $P_1$-Anosov representation is a quasi-isometric embedding.

The action of a non-elementary Gromov hyperbolic group on its boundary is minimal, and the attracting fixed points of its infinite-order elements are dense \cite[Proposition~4.2]{KB02}. We shall use these facts repeatedly. From now on we write $\xi=\xi_\rho$ and $\Lambda=\Lambda_\rho$.

\subsection{Negative spheres and the graph model}\label{sec:graph}

Let $p>q\ge2$, let $\Gamma$ be a non-elementary Gromov hyperbolic group with $\partial_\infty\Gamma\cong\mathbb S^{p-1}$, and let $\rho:\Gamma\longrightarrow\PO_{p,q}$ be $\HH^{p,q-1}$-convex cocompact. Write $\Lambda=\Lambda_\rho$.

Transversality makes the boundary map injective, and hence it is a topological embedding because $\partial_\infty\Gamma$ is compact. Thus $\Lambda$ is homeomorphic to $\mathbb S^{p-1}$. Since $p>q\ge2$, this sphere is simply connected. The restriction of the positive-ray double cover to $\Lambda$ is therefore trivial, and its full preimage has two components exchanged by the antipodal map. Choose a lift of $\Lambda$ with negative pairings. This lift is continuous: locally, its sign is determined by negative pairing with a fixed lift of another point of $\Lambda$. Its image is therefore one of the two components, and its antipode has the same negative-pairing property.

Using the sphere-product model $\partial_\infty\widehat{\HH}^{p,q-1}\cong\mathbb S^{p-1}\times\mathbb S^{q-1}$, Lemma~3.3 and Corollary~3.4 of \cite{BK25} show that either component is the graph of a strictly $1$-Lipschitz map $\varphi:\mathbb S^{p-1}\longrightarrow\mathbb S^{q-1}$. Here ``strictly'' means that
$$
d_{\mathbb S^{q-1}}(\varphi(x),\varphi(y))
<d_{\mathbb S^{p-1}}(x,y)
\quad\text{for all }x\ne y;
$$

it does not mean that the Lipschitz constant is uniformly smaller than $1$.

\begin{lemma}\label{lem:c1graph}
If $\Lambda$ is a $C^1$ submanifold, then $\varphi$ is $C^1$.
\end{lemma}

\begin{proof}
Work on one lifted component $\widetilde\Lambda\subset\mathbb S^{p-1}\times\mathbb S^{q-1}$. It is a $C^1$ submanifold, and the graph projection $\pi_+:\widetilde\Lambda\to\mathbb S^{p-1}$ is a homeomorphism. We first claim that the boundary conformal form is nonnegative on $T\widetilde\Lambda$. Let $c(t)$ be a $C^1$ curve in $\widetilde\Lambda$, represented by normalized null vectors, and put $z=c(0)$. For $t\ne0$ sufficiently small, negativity gives $\langle c(t),z\rangle\le0$, with strict inequality whenever $c(t)\ne z$. Since both $c(t)$ and $z$ are null,
$$
\left\langle\frac{c(t)-z}{t},\frac{c(t)-z}{t}\right\rangle
=-\frac{2}{t^2}\langle c(t),z\rangle\ge0.
$$
Passing to the limit shows that the boundary conformal form is nonnegative on $T\widetilde\Lambda$.

A nonzero vertical vector has negative norm for $g_{\mathbb S^{p-1}}-g_{\mathbb S^{q-1}}$, so $T\widetilde\Lambda$ contains no such vector. Hence $d\pi_+$ is injective at every point. The source and target have the same dimension, so $d\pi_+$ is an isomorphism everywhere. The inverse function theorem makes $\pi_+$ a local $C^1$ diffeomorphism. Since it is globally a homeomorphism, its inverse is $C^1$, and $\varphi=\pi_-\circ\pi_+^{-1}$ is $C^1$.
\end{proof}

\subsection{The tangent form and its radical}\label{sec:tangent}
Suppose that $\Lambda$ is a $C^1$ submanifold. By Lemma~\ref{lem:c1graph}, the map $\varphi$ is $C^1$. At $x\in\mathbb S^{p-1}$ define
$$
Q_x(v,w)=\ip{v}{w}_{\mathbb S^{p-1}}
-\ip{d\varphi_xv}{d\varphi_xw}_{\mathbb S^{q-1}}.
$$

Under the graph parametrization $x\mapsto(x,\varphi(x))$, the form $Q$ is the pullback of $g_{\mathbb S^{p-1}}-g_{\mathbb S^{q-1}}$. Put $A_x=\id-(d\varphi_x)^*d\varphi_x$. Differentiating the $1$-Lipschitz inequality along geodesics gives $\lVert d\varphi_xv\rVert\le\lVert v\rVert$, so $A_x$ is positive semidefinite. Moreover, $Q_x(v,w)=\langle A_xv,w\rangle_{\mathbb S^{p-1}}$. The radical of $Q_x$ is
$$
R_x=\Rad(Q_x)=\{v\in T_x\mathbb S^{p-1}:(d\varphi_x)^*d\varphi_xv=v\}.
$$

Thus $R_x=\ker A_x$ is the eigenspace of $(d\varphi_x)^*d\varphi_x$ corresponding to the eigenvalue $1$, and $d\varphi_x|_{R_x}:R_x\longrightarrow T_{\varphi(x)}\mathbb S^{q-1}$ is an isometric injection. Therefore $\dim R_x\le q-1$.

Glorieux--Monclair use this form in the AdS case to find a spacelike tangent point and then propagate spacelikeness by minimality \cite[Lemma~4.1 and Corollary~4.2]{GM24}; immediately afterward they explain why the first step fails in higher signature. Monclair also observes that the restricted conformal form on a maximal-dimensional $C^1$ limit set has constant signature and describes a resulting real-rank bound \cite[paragraph following Theorem~3.11, p.~37]{Mon25}. Here we study the radical subbundle itself: its topology and its behavior under the derivative action allow us to treat tangent spaces that are not initially known to be spacelike. This pointwise radical should be distinguished from the radical of the ambient form on $\operatorname{span}(\Lambda)$, considered in \cite[Definition~2.8, Lemmas~3.7 and~3.9, and Proposition~3.10]{BK25}.

\begin{proposition}\label{prop:constant}

In the preceding general-signature setting, the dimension of $R_x$ is independent of $x$. Moreover, $R:=\bigsqcup_{x\in\mathbb S^{p-1}}R_x$ is a continuous vector subbundle of $T\mathbb S^{p-1}$. Its rank is therefore a well-defined integer $r\in\{0,\ldots,q-1\}$. In particular, when $q=3$ one has $r\in\{0,1,2\}$.

\end{proposition}

\begin{proof}

Let $g_{\operatorname{Ein}}$ denote the representative of the Einstein conformal class determined by the sphere-product model above. For $z\in\Lambda$, let $\mathcal R_z:=\Rad\left(g_{\operatorname{Ein}}|_{T_z\Lambda}\right)$ and put $\nu(z)=\dim\mathcal R_z$. This does not depend on the representative chosen from the Einstein conformal class. Under the graph parametrization, $\mathcal R_z$ is identified with $R_x$. It is therefore enough to prove that $\nu$ is constant.

The action of $\PO_{p,q}$ on the Einstein universe is conformal. Thus, for every $\gamma\in\Gamma$ and $z\in\Lambda$, there is a number $c_{\gamma,z}>0$ such that
$$
g_{\operatorname{Ein},\,\rho(\gamma)z}
\left(d(\rho(\gamma))_z u,d(\rho(\gamma))_z v\right)
=c_{\gamma,z}\,g_{\operatorname{Ein},\,z}(u,v)
$$
for all $u,v\in T_z\partial_\infty\HH^{p,q-1}$. Since $\rho(\gamma)$ preserves $\Lambda$, its derivative maps $T_z\Lambda$ onto $T_{\rho(\gamma)z}\Lambda$. It follows that $d(\rho(\gamma))_z(\mathcal R_z)=\mathcal R_{\rho(\gamma)z}$. In particular, $\nu$ is $\rho(\Gamma)$-invariant.

For $k\ge0$, set $\mathcal X_k:=\{z\in\Lambda:\nu(z)\ge k\}$. We claim that $\mathcal X_k$ is closed. Indeed, suppose that $z_j\in\mathcal X_k$ and $z_j\to z$. Choose a local trivialization of $T\Lambda$ near $z$ and a local representative of the Einstein conformal class. For each $j$, choose a $k$-plane $E_j\subset\mathcal R_{z_j}$. After passing to a subsequence, the compactness of the Grassmannian gives a limiting $k$-plane $E\subset T_z\Lambda$. The tangent spaces of the $C^1$ submanifold $\Lambda$ vary continuously, as do the restrictions of the conformal form to them. Hence $g_{\operatorname{Ein},\,z}(u,v)=0$ for every $u\in E$ and $v\in T_z\Lambda$. Thus $E\subset\mathcal R_z$, and $z\in\mathcal X_k$.

The boundary map is an equivariant homeomorphism from $\partial_\infty\Gamma$ onto $\Lambda$. Since the action of a non-elementary Gromov hyperbolic group on its boundary is minimal, the action of $\rho(\Gamma)$ on $\Lambda$ is minimal as well. Each $\mathcal X_k$ is closed and invariant, and therefore is either empty or equal to $\Lambda$.

Fix $z_0\in\Lambda$ and let $r=\nu(z_0)$. Since $z_0\in\mathcal X_r$, minimality gives $\mathcal X_r=\Lambda$. On the other hand, $z_0\notin\mathcal X_{r+1}$, so $\mathcal X_{r+1}$ cannot be all of $\Lambda$ and must be empty. It follows that $\nu(z)=r$ for every $z\in\Lambda$. Equivalently, $\dim R_x=r$ for every $x\in\mathbb S^{p-1}$. The estimate above gives $r\le q-1$.

It remains to check that the spaces $R_x$ form a continuous vector bundle. Since $\varphi$ is $C^1$, the endomorphism $A_x=\id-(d\varphi_x)^*d\varphi_x$ depends continuously on $x$. We have $R_x=\ker A_x$, and $\rank A_x=p-1-r$ is constant.

Fix $x_0\in\mathbb S^{p-1}$ and choose a local frame for $T\mathbb S^{p-1}$ near $x_0$. In this frame, $A_x$ is a continuous family of matrices of rank $p-1-r$. Some $(p-1-r)\times(p-1-r)$ minor of $A_{x_0}$ is nonzero, and it remains nonzero in a neighborhood $U$ of $x_0$. After reordering rows and columns, the corresponding independent equations have the form $B(x)u+C(x)v=0$, where $B(x)$ is invertible and $v\in\R^r$. Since the full matrix has rank $p-1-r$, these equations determine its kernel. Thus,
$$
\ker A_x=\left\{\bigl(-B(x)^{-1}C(x)v,v\bigr):v\in\R^r\right\}
$$
for $x\in U$. Taking the standard basis vectors for $v$ gives a continuous local frame for $\ker A$. Hence $R=\ker A$ is a continuous rank-$r$ vector subbundle of $T\mathbb S^{p-1}$.
\end{proof}

\begin{remark}

The $C^1$ assumption enters the proof through the continuity of $d\varphi$. If the graph is only Lipschitz, then $d\varphi$ is defined almost everywhere, and the preceding minimality argument does not give a continuous radical distribution.

\end{remark}

\section{The extreme radical ranks}\label{sec:r0r2}

\subsection{The spacelike case \texorpdfstring{$r=0$}{r=0}}

When $r=0$, the limit set $\Lambda$ is spacelike, so rigidity follows from Monclair's theorem \cite[Theorem~3.11]{Mon25}. For completeness, we recall a proof using the critical-exponent formula of Glorieux--Monclair \cite[Theorem~5.5]{GM21} and the entropy rigidity theorem of Mazzoli--Viaggi \cite[Theorem~4]{MV26}.

We recall the pseudo-Riemannian critical exponent used below. Put $\Delta=\rho(\Gamma)$, and let $d_{\HH^{p,q-1}}$ be the invariant pseudo-distance which is the length of the spacelike geodesic segment joining two space-related points, and is zero otherwise. For a basepoint $o$ in the convex hull of $\Lambda$, the pseudo-Riemannian critical exponent of $\Delta$ is
$$
\delta_\Delta
=\limsup_{R\longrightarrow\infty}\frac1R
\log\#\{\gamma\in\Delta:
d_{\HH^{p,q-1}}(o,\gamma o)\le R\}.
$$

It is independent of the choice of $o$; see \cite[Theorem~3.5 and Definition~3.6]{GM21}.

\begin{proposition}\label{prop:r0entropy} 
In the general-signature setting above, if $r=0$, then $\rho(\Gamma)$ preserves a totally geodesic copy of $\HH^p$.
\end{proposition}
\begin{proof}
Since $r=0$, the restriction of the boundary conformal form to $T\Lambda$ is positive definite. Compactness makes it uniformly positive definite with respect to any auxiliary Riemannian metric on $\Lambda$.

Fix a point $o$ in the convex hull of $\Lambda$, and let $d_o$ be the boundary Gromov quasi-distance of Glorieux--Monclair \cite{GM21}. Choose a unit lift $\widehat o$ of $o$, and for each $\xi\in\Lambda$ choose the unique lift $\widehat\xi$ satisfying $\langle\widehat\xi,\widehat o\rangle=-1$. In this normalization, the formula in \cite[Section~3.5]{GM21} becomes
$$
d_o(\xi,\eta)=\frac12 \sqrt{\langle\widehat\xi-\widehat\eta, \widehat\xi-\widehat\eta\rangle}.
$$

Choose a $C^1$ local parametrization $u\mapsto z(u)$ of $\Lambda$, and write $F(u)=\widehat{z(u)}$ for its normalized lift. The quadratic form $\langle dF_u(\cdot),dF_u(\cdot)\rangle$ represents the restricted boundary conformal structure and is positive definite. On a sufficiently small relatively compact coordinate neighborhood, it is uniformly positive definite, and $C^1$ regularity gives $F(v)-F(u)=dF_u(v-u)+o(|v-u|)$ uniformly as $|v-u|\to0$. Substituting into the formula for $d_o$ gives constants $c,C>0$ such that $c|u-v|\le d_o(z(u),z(v))\le C|u-v|$ for all sufficiently close $u,v$. Thus $d_o$ is locally bi-Lipschitz equivalent on $\Lambda$ to an auxiliary Riemannian distance, and $\operatorname{Hdim}_{d_o}(\Lambda)=p-1$. By \cite[Theorem~5.5]{GM21}, $\delta_\Delta=p-1$.

By Definition~\ref{def:hcc}, $K:=\ker\rho$ is finite. Hence $\Delta\cong\Gamma/K$, and the quotient map, being a quasi-isometry, induces a homeomorphism $\partial_\infty\Gamma\cong\partial_\infty\Delta$. Let $\widetilde\Delta<\Ogrp_{p,q}$ be the full inverse image of $\Delta$. This group is finitely generated and linear, so Selberg's lemma \cite{Sel60} gives a torsion-free subgroup $\widetilde\Delta_0<\widetilde\Delta\cap\SO_{p,q}$ of finite index in $\widetilde\Delta$. Since the kernel of projectivization is finite, projectivization is injective on $\widetilde\Delta_0$. Its image $\Delta_0<\Delta$ therefore has finite index and admits a lift to $\SO_{p,q}$. This subgroup remains $\HH^{p,q-1}$-convex cocompact and has the same Gromov boundary, proximal limit set, and critical exponent.

By \cite[Corollary~1.1]{SST}, the group $\Delta_0$ preserves a unique connected complete maximal spacelike $p$-submanifold $M\subset\HH^{p,q-1}$ with $\partial_\infty M=\Lambda$, and acts properly discontinuously and cocompactly on it. By \cite[Lemma~3.11]{SST}, the submanifold $M$ is an entire spacelike $p$-graph, hence is diffeomorphic to $\R^p$. Since $\Delta_0$ is torsion-free, its properly discontinuous action on $M$ is free: proper discontinuity makes every point stabilizer finite. After replacing $\Delta_0$ by a subgroup of index at most two, we may therefore assume that $N:=M/\Delta_0$ is a closed orientable aspherical $p$-manifold with $\pi_1(N)\cong\Delta_0$.

Since the critical exponent is unchanged under passage to a finite-index subgroup, $\delta_{\Delta_0}=p-1$. The lift of $\Delta_0$ to $\SO_{p,q}$, viewed as a representation of $\pi_1(N)$, therefore satisfies the equality case of \cite[Theorem~4]{MV26}. Hence $\Delta_0$ preserves a totally geodesic subspace $X\subset\HH^{p,q-1}$ isometric to $\HH^p$. Since $X$ is itself a complete maximal spacelike $p$-submanifold, uniqueness in \cite[Corollary~1.1]{SST} gives $X=M$.

Finally, $\partial_\infty X=\Lambda$, and $\Delta$ preserves $\Lambda$. Since a totally geodesic copy of $\HH^p$ is determined by its ideal boundary, the whole group $\Delta=\rho(\Gamma)$ preserves $X$.
\end{proof}

\subsection{The maximal radical rank}

\begin{proposition}\label{prop:maxradpq}
Let $p>q\ge2$. In the general-signature setting of Subsections~\ref{sec:graph} and~\ref{sec:tangent}, the case $r=q-1$ cannot occur.
\end{proposition}

\begin{proof}
Suppose that $r=q-1$. The restriction $d\varphi_x|_{R_x}:R_x\longrightarrow T_{\varphi(x)}\mathbb S^{q-1}$ is an isometric injection between spaces of the same dimension, so $d\varphi_x$ is surjective for every $x$. Thus $\varphi$ is a $C^1$ submersion. Its image is open, and compactness of the domain makes it closed. Since $\mathbb S^{q-1}$ is connected, $\varphi$ is surjective.

Choose $x,z\in\mathbb S^{p-1}$ such that $\varphi(z)=-\varphi(x)$. The strict distance inequality from Subsection~\ref{sec:graph} gives
$$
\pi=d_{\mathbb S^{q-1}}(\varphi(x),\varphi(z))<d_{\mathbb S^{p-1}}(x,z)\le\pi,
$$
a contradiction.
\end{proof}

We retain the general-signature setting in the next section. By Proposition~\ref{prop:maxradpq}, the radical rank satisfies $0\le r\le q-2$. When $q=3$, this leaves $r\in\{0,1\}$. The proof of Theorem~\ref{thm:main} below does not use Proposition~\ref{prop:r0entropy}, since the same algebraic argument treats both possibilities.

\section{Spectral faces and the Zariski closure}\label{sec:wall}

The representation $\rho:\Gamma\to\PO_{p,q}$ need not lift to $\SO_{p,q}$. We therefore let $\widetilde\Delta<\Ogrp_{p,q}$ be the full inverse image of $\rho(\Gamma)$ and take its real Zariski closure in $\Ogrp_{p,q}$. Multiplication by $-\id$ changes neither singular values nor eigenvalue moduli, so this finite central extension does not change the Cartan and Jordan coordinates, the Benoist limit cone, the limit set, or invariant projective subspaces.

For a subgroup $\Delta<\Ogrp_{p,q}$, its \emph{real Zariski closure} is the smallest real algebraic subgroup containing $\Delta$; we denote it by $\Zcl(\Delta)$. If $H$ is a real algebraic group, $H^\circ$ denotes its algebraic identity component. For a projective subgroup, $\Zcl$ denotes the real Zariski closure of its full inverse image in the corresponding orthogonal group. Since the determinant has finite image, the algebraic identity component of every real algebraic subgroup of $\Ogrp_{p,q}$ is contained in $\SO_{p,q}$. Thus $\Ogrp_{p,q}$ is used only to accommodate the possibly disconnected lift, whereas the ambient linear group for the connected Zariski-closure and restricted-weight arguments below is $G=\SO_{p,q}$.

We use the chamber, coordinates, and faces of $G$ introduced in Subsection~\ref{sec:projections}; they agree with the folded coordinates for $\PO_{p,q}$ and $\Ogrp_{p,q}$. We write the coordinates of $\lambda_G(g)$ as $\lambda_1(g)\ge\cdots\ge\lambda_q(g)\ge0$ and, for $1\le m\le q-1$, set
$$
F_m=\mathfrak a_G^+\cap\{(t_1,\ldots,t_q)\in\mathbb R^q:t_{m+1}=0\}.
$$

Equivalently, $\lambda_G(g)\in F_m$ if and only if $\lambda_{m+1}(g)=0$.

\subsection{The face forced by the tangent radical}

\begin{lemma}\label{lem:spectral}
Let $\gamma\in\Gamma$ have infinite order and let $g=\rho(\gamma)$. Then
$$
\#\{i\in\{2,\ldots,q\}:\lambda_i(g)>0\}\le r.
$$
In particular, $\lambda_{r+2}(g)=0$.
\end{lemma}

\begin{proof}
Let $\ell=\xi(\gamma^+)$ be the attracting isotropic line. Choose a lift $\widetilde g\in\Ogrp_{p,q}$ of $g$ and a nonzero $x_+\in\ell$ such that $\widetilde g x_+=a x_+$ for some $a\in\mathbb R^*$. Let $\bar g$ denote the action induced by $\widetilde g$ on $W_\ell:=\ell^\perp/\ell$. It preserves the form of signature $(p-1,q-1)$ and has eigenvalue moduli
$$
e^{\lambda_2(g)},\ldots,e^{\lambda_q(g)},\ 
\underbrace{1,\ldots,1}_{p-q},\ 
e^{-\lambda_q(g)},\ldots,e^{-\lambda_2(g)}.
$$

Under the standard identification $T_\ell\partial_\infty\HH^{p,q-1}\cong\operatorname{Hom}(\ell,W_\ell)$, evaluation at $x_+$ identifies this tangent space with $W_\ell$. In these coordinates the projective derivative is $a^{-1}\bar g$. We therefore regard $E_\ell=T_\ell\Lambda$ as a $(p-1)$-dimensional nonnegative subspace of $W_\ell$. Its invariance under the derivative is equivalent to its invariance under $\bar g$. Write $\mathcal R_\ell$ for the radical of the form restricted to $E_\ell$.

The map $\bar g$ preserves the invariant filtration
$$
0\subset\mathcal R_\ell\subset\mathcal R_\ell^\perp\subset W_\ell.
$$

Since $\mathcal R_\ell$ is isotropic of dimension $r$, the quotient $\mathcal R_\ell^\perp/\mathcal R_\ell$ has signature $(p-1-r,q-1-r)$. The subspace $E_\ell/\mathcal R_\ell$ is positive definite of dimension $p-1-r$, so its orthogonal complement $E_\ell^\perp/\mathcal R_\ell$ is negative definite of dimension $q-1-r$. Thus $E_\ell^\perp$ is nonpositive with radical $\mathcal R_\ell$, and
$$
\mathcal R_\ell^\perp/\mathcal R_\ell
= (E_\ell/\mathcal R_\ell)
\mathbin{\overset{\perp}{\oplus}}(E_\ell^\perp/\mathcal R_\ell).
$$

Thus every eigenvalue of $\bar g$ on $\mathcal R_\ell^\perp/\mathcal R_\ell$ has modulus $1$. The ambient form induces a perfect pairing between $\mathcal R_\ell$ and $W_\ell/\mathcal R_\ell^\perp$, so the actions on these spaces are contragredient. Their eigenvalues can therefore be listed in $r$ reciprocal pairs, counted with algebraic multiplicity. Each pair contributes at most one eigenvalue of modulus greater than $1$. Hence $\bar g$ has at most $r$ such eigenvalues. Comparing with the list above proves the result.
\end{proof}

\begin{corollary}\label{cor:wall}

The Benoist limit cone of $\rho(\Gamma)$ is contained in $F_{r+1}$. In particular, it is contained in the last Weyl wall $F_{q-1}$, and $\rho(\Gamma)$ is not Zariski dense in $\PO_{p,q}$.

\end{corollary}

\begin{proof}

The Jordan projection of a finite-order element is zero, so the first statement follows from Lemma~\ref{lem:spectral}. Since $r\le q-2$, the face $F_{r+1}$ has empty interior in $\mathfrak a_G^+$. If $\rho(\Gamma)$ were Zariski dense, its intersection with the algebraic identity component would have finite index and be Zariski dense there. If $N$ is the exponent of the finite component group, then $\delta^N$ belongs to this intersection and its Jordan projection is $N$ times that of $\delta$ for every $\delta\in\rho(\Gamma)$; hence passing to the identity component does not change the limit cone. Benoist's theorem \cite[Theorem~1.2]{Ben97} would imply that this cone has nonempty interior, a contradiction.
\end{proof}

\subsection{Consequences for the Zariski closure}

Zariski density does not imply equality between the limit cone of a subgroup and that of its Zariski closure. Indeed, for a connected reductive group $L$ one has $\limcone_L^L=\mathfrak a_L^+$, whereas the intrinsic cone $\limcone_\Delta^L$ of a Zariski-dense subgroup $\Delta<L$ need only have nonempty interior and may be a proper full-dimensional subcone of $\mathfrak a_L^+$ \cite[Theorem~1.2(b)]{Ben97}. The following lemma shows that containment in one of the special faces $F_m$ nevertheless passes from $\Delta$ to $L$.

\begin{lemma}\label{lem:wall-inheritance}

Let $G=\SO_{p,q}$ and let $1\le m\le q-1$. Let $L<G$ be a connected semisimple real algebraic subgroup, and let $\Delta<L$ be Zariski dense. Assume that $\limcone_\Delta^G\subset F_m$. Let $A_L<L$ be a maximal real split torus and put $\mathfrak a_L=\operatorname{Lie}(A_L)$. For the standard $L$-module $V=\mathbb R^{p,q}$, set $V_\chi:=\{v\in V:Xv=\chi(X)v\text{ for all }X\in\mathfrak a_L\}$ and $\Omega_L(V):=\{\chi\in\mathfrak a_L^*:V_\chi\ne0\}$. Thus
$$
V=\bigoplus_{\chi\in\Omega_L(V)}V_\chi.
$$

Then, for every $X\in\mathfrak a_L$,
$$
\sum_{\substack{\chi\in\Omega_L(V)\\\chi(X)>0}}\dim V_\chi\le m
$$
and $\rank_\mathbb R L\le m$. Moreover, $\limcone_L^G\subset F_m$.
\end{lemma}

\begin{proof}
Benoist's theorem \cite[Theorem~1.2]{Ben97} says that the intrinsic limit cone $\limcone_\Delta^L$ has nonempty interior in $\mathfrak a_L^+$. After conjugating $A_L$ into a maximal split torus of $G$, define the continuous, positively homogeneous folding map
$$
\iota_*:\mathfrak a_L^+\longrightarrow\mathfrak a_G^+,
\qquad X\longmapsto\lambda_G(\exp X).
$$

For every $\delta\in\Delta$, functoriality of the multiplicative Jordan decomposition gives $\iota_*(\lambda_L(\delta))=\lambda_G(\delta)$. Hence $\iota_*$ sends $\limcone_\Delta^L$ into $\limcone_\Delta^G\subset F_m$. Choose $A$ in the interior of $\limcone_\Delta^L$ and outside the finitely many hyperplanes $\ker\chi$, where $\chi$ ranges over the nonzero weights in $\Omega_L(V)$. Since $\iota_*(A)\in F_m$, at most $m$ logarithms of eigenvalue moduli of $\exp A$ are positive. Equivalently,
$$
\sum_{\substack{\chi\in\Omega_L(V)\\\chi(A)>0}}\dim V_\chi\le m.
$$

The orthogonal form pairs $V_\chi$ nondegenerately with $V_{-\chi}$, so opposite weights have equal multiplicity. Because $A$ is \emph{weight-regular}, meaning that $\chi(A)\ne0$ for every nonzero restricted weight $\chi\in\Omega_L(V)$, exactly one weight in each opposite pair is positive on $A$. Thus the sum above equals $\frac12\sum_{\chi\ne0}\dim V_\chi$. For arbitrary $X\in\mathfrak a_L$, at most one weight in each opposite pair is positive on $X$, so the same upper bound $m$ holds. Since the inclusion $L<\SO(V)$ is faithful, the nonzero restricted weights span $\mathfrak a_L^*$. There are at most $m$ opposite pairs, hence $\rank_\mathbb R L=\dim\mathfrak a_L\le m$.

For every $g\in L$, its hyperbolic part is $L$-conjugate to $\exp X$ for some $X\in\mathfrak a_L$. The ambient eigenvalue moduli are determined by the numbers $\chi(X)$, so the preceding inequality gives $\lambda_{m+1}(g)=0$. Thus every ray used to define $\limcone_L^G$ lies in $F_m$, and therefore $\limcone_L^G\subset F_m$.
\end{proof}

\begin{lemma}\label{lem:acausal-dim}
Let $U$ be a real vector space of signature $(a,b)$ with $a,b\ge1$, and put
$$
\partial_\infty\HH(U)=\{[u]\in\Pj(U):\langle u,u\rangle=0\}.
$$
Any compact negative subset of $\partial_\infty\HH(U)$ has topological dimension at most $a-1$.
\end{lemma}

\begin{proof}
The assertion is immediate if the set has at most one point. Choose an orthogonal splitting $U=U_+\oplus U_-$ with $\dim U_+=a$ and $\dim U_-=b$. Choose a negative lift $\widetilde A$ of the given set to the positive-ray cover and normalize it in $\mathbb S(U_+)\times\mathbb S(U_-)$. This lift is continuous: locally, its sign is determined by negative pairing with a fixed lift of another point of the set. Hence $\widetilde A$ is compact and projects homeomorphically to the original set. The projection $\widetilde A\to\mathbb S(U_+)$ is injective. Indeed, two normalized points $(u,w_1)$ and $(u,w_2)$ with the same positive coordinate have pairing $1-\langle w_1,w_2\rangle\ge0$, which contradicts negativity unless they coincide. The continuous injection $\widetilde A\to\mathbb S(U_+)$ is therefore a topological embedding. Since $\dim\mathbb S(U_+)=a-1$, the result follows.
\end{proof}

\begin{lemma}\label{lem:noiso}
Let $p>q\ge2$ and let $H=\Zcl(\rho(\Gamma))<\Ogrp_{p,q}$. Then $H$ preserves no nonzero totally isotropic subspace $I<V=\R^{p,q}$.
\end{lemma}

\begin{proof}
Assume that $I$ is nonzero, totally isotropic, and $H$-invariant. Then $I^\perp$ is also invariant, so $A:=\Lambda\cap\Pj(I^\perp)$ is closed and $\Gamma$-invariant. By minimality, $A=\varnothing$ or $A=\Lambda$.

If $A=\Lambda$, first put $B=\Lambda\cap\Pj(I)$. This is again closed and $\Gamma$-invariant, hence empty or all of $\Lambda$. The latter is impossible: $\Lambda$ contains at least two points, whereas any two vectors in the totally isotropic space $I$ pair to zero. Thus $B=\varnothing$.

Put $d=\dim I$. Projection induces a continuous map $\pi:\Lambda\longrightarrow\partial_\infty\HH(I^\perp/I)$. It is injective. Indeed, if $\pi([x])=\pi([y])$, then, after rescaling, $x-cy\in I$ for some $c\ne0$. Since $x,y\in I^\perp$ and both are isotropic,
$$
0=\langle x-cy,x-cy\rangle=-2c\langle x,y\rangle.
$$

Thus, $\langle x,y\rangle=0$, and negativity forces $[x]=[y]$. If $d=q$, then $I^\perp/I$ is positive definite, contradicting the existence of the nonzero isotropic line $\pi([x])$. If $d<q$, the quotient has signature $(p-d,q-d)$. Its induced form preserves the pairings of vectors in $I^\perp$, so $\pi(\Lambda)$ is negative. Since $\Lambda$ is compact, $\pi$ is a topological embedding, and Lemma~\ref{lem:acausal-dim} gives
$$
p-1=\dim\Lambda\le p-d-1,
$$
again a contradiction.

If $A=\varnothing$, choose an arbitrary infinite-order element $\gamma\in\Gamma$, put $g=\rho(\gamma)$, and let $\ell^\pm=\xi(\gamma^\pm)$. By the $P_1$-Anosov property, $g$ and $g^{-1}$ are proximal, with attracting lines $\ell^+$ and $\ell^-$, respectively. Since $\ell^+\not\subset I^\perp$, some $v\in I$ pairs nontrivially with a lift of $\ell^+$. For an orthogonal proximal element, the repelling hyperplane of $g^{-1}$ is $\Pj((\ell^+)^\perp)$. Thus $[v]$ lies outside this hyperplane, and projective proximal dynamics gives $g^{-n}[v]\to\ell^-$. Since $\Pj(I)$ is invariant and closed, $\ell^-\subset I\subset I^\perp$. But $\ell^-\in\Lambda$, contradicting $A=\varnothing$.
\end{proof}

\begin{corollary}\label{cor:red}
Under the same general-signature hypotheses, let $H=\Zcl(\rho(\Gamma))<\Ogrp_{p,q}$. Then the algebraic identity component $H^\circ$ is reductive. Moreover, every $H$-invariant subspace of $V=\R^{p,q}$ is nondegenerate.
\end{corollary}

\begin{proof}
Let $U=R_u(H^\circ)$ be the unipotent radical of $H^\circ$. Since $H^\circ$ is normal in $H$ and $U$ is characteristic in $H^\circ$, the group $H$ normalizes $U$. Hence $F:=V^U$ is $H$-invariant. Every element of $\mathfrak u:=\operatorname{Lie}(U)$ acts nilpotently on $V$. Engel's theorem \cite[Theorem~1.35]{Knapp} therefore gives a nonzero vector annihilated by $\mathfrak u$. Since $U$ is connected, this vector is fixed by $U$, and hence $F\ne0$.

Suppose first that $F$ is degenerate. Then $\Rad(F):=F\cap F^\perp$ is a nonzero $H$-invariant totally isotropic subspace, contrary to Lemma~\ref{lem:noiso}. Thus $F$ is nondegenerate and $V=F\oplus F^\perp$.

The subspace $F^\perp$ is $U$-invariant. If $F^\perp\ne0$, the same argument gives a nonzero $U$-fixed vector in $F^\perp$. This is impossible, because every $U$-fixed vector belongs to $F$, while $F\cap F^\perp=0$. Hence $F^\perp=0$ and $F=V$. The group $U$ therefore acts trivially on $V$, and the faithful inclusion $H<\Ogrp_{p,q}$ gives $U=\{1\}$. Thus $H^\circ$ is reductive.

Finally, were an $H$-invariant subspace $W<V$ degenerate, then $W\cap W^\perp$ would~be a nonzero $H$-invariant totally isotropic subspace, contradicting Lemma~\ref{lem:noiso}.
\end{proof}

\begin{remark}
Corollary~\ref{cor:red} establishes reductivity but does not exclude a positive-dimensional center. In the proof of Theorem~\ref{thm:main}, we will also need semisimplicity when the Zariski closure acts irreducibly. Since the negative index is three, this follows from the next lemma.
\end{remark}

We shall use {\em Schur's lemma} in the following standard form. If $V$ and $W$ are irreducible modules, every equivariant map $V\to W$ is either zero or an isomorphism; in particular, $\operatorname{End}_H(V)$ is a division algebra. Over $\mathbb C$ this algebra consists of scalars, while for an irreducible real module it is isomorphic to $\mathbb R$, $\mathbb C$, or $\mathbb H$. We refer to \cite[Proposition~5.1 and Corollary~5.2]{Knapp} for Schur's lemma and to \cite[Section~2.1]{DSL11} for the real, complex, and quaternionic alternatives used below.

\begin{lemma}[Semisimplicity in odd index]\label{lem:semisimple-irreducible}
Let $L<\Ogrp_{a,b}$ be connected and reductive, acting irreducibly on $V=\R^{a,b}$. If at least one of $a$ and $b$ is odd, then the Lie algebra of $L$ is semisimple.
\end{lemma}

\begin{proof}
Suppose that the center of $\operatorname{Lie}(L)$ contains a nonzero element $Z$. Since $L$ is connected, $Z$ commutes with its action on $V$, so $Z\in D:=\operatorname{End}_L(V)$. Schur's lemma makes $D$ a division algebra. The minimal polynomial of $Z$ over $\R$ is therefore irreducible: a nontrivial factorization would give a product of two nonzero elements of $D$ equal to zero. Since $Z$ is skew-adjoint for the invariant symmetric form, $\operatorname{tr}Z=0$. It cannot be a nonzero real scalar, so its minimal polynomial is an irreducible quadratic. Its complex roots have the form $\alpha\pm i\beta$ with $\beta>0$; their multiplicities agree because $Z$ is real, and the trace condition gives $\alpha=0$. Thus $Z^2=-\beta^2\id$.

Set $J=\beta^{-1}Z$. Then $J^2=-\id$ and $J^*=-J$, so $J$ is an orthogonal complex structure. The invariant symmetric form is therefore the real part of a nondegenerate Hermitian form, so both its positive and negative indices are even, contrary to the hypothesis. Hence the Lie algebra of the center is zero, and reductivity implies that $\operatorname{Lie}(L)$ is semisimple.
\end{proof}

\section{The reducible case}\label{sec:descent}

\begin{lemma}\label{lem:directsum}
Let $\Gamma$ be a non-elementary Gromov hyperbolic group with connected boundary, and let $\sigma:\Gamma\to\PGL(U)$ be a representation. Suppose that there is a continuous map $\xi:\partial_\infty\Gamma\to\Pj(U)$ such that, for every infinite-order $\gamma\in\Gamma$, the element $\sigma(\gamma)$ is proximal with attracting line $\xi(\gamma^+)$. Suppose also that $U=U_1\oplus\cdots\oplus U_k$ and that every element of $\sigma(\Gamma)$ preserves each projective subspace $\Pj(U_i)$. Then $\xi(\partial_\infty\Gamma)$ is contained in one of the subspaces $\Pj(U_i)$.
\end{lemma}

\begin{proof}
For every infinite-order $\gamma\in\Gamma$, the simple attracting eigenline of $\sigma(\gamma)$ lies in one of the invariant blocks. The points $\gamma^+$ are dense in $\partial_\infty\Gamma$, and $\xi$ is continuous. Hence $\xi(\partial_\infty\Gamma)$ is contained in the finite union $\bigcup_i\Pj(U_i)$. Since $\partial_\infty\Gamma$ is connected and these projective subspaces are closed and disjoint, the image lies in a single block.
\end{proof}

\begin{proposition}\label{prop:descent}
Let $p>q\ge2$, and let $\rho$ satisfy the general-signature hypotheses of Subsections~\ref{sec:graph} and~\ref{sec:tangent}. Let $H=\Zcl(\rho(\Gamma))<\Ogrp_{p,q}$. If $H$ is reducible on $V=\R^{p,q}$, then $\rho(\Gamma)$ preserves a nondegenerate proper subspace $W<V$ of signature $(p,b)$, where $1\le b\le q-1$, and $W^\perp$ is negative definite.

Restriction to $W$ induces a representation $\rho_W:\Gamma\longrightarrow\PO(W)\cong\PO_{p,b}$, which is $\HH^{p,b-1}$-convex cocompact and has limit set $\Lambda\subset\Pj(W)$. This limit set is again a $C^1$ submanifold, and its tangent form and radical rank are unchanged. Moreover, $\Zcl(\rho_W(\Gamma))$ acts irreducibly on $W$.
\end{proposition}
\begin{proof}
By Corollary~\ref{cor:red}, every $H$-invariant subspace is nondegenerate. Since $H$ is reducible, there is an $H$-invariant orthogonal decomposition $V=U\oplus U^\perp$ with both summands nonzero. Lemma~\ref{lem:directsum} implies that $\Lambda$ lies in the projectivization of one of these two summands. Therefore $W:=\operatorname{span}_{\R}\{v\in V\setminus\{0\}:[v]\in\Lambda\}$ is a proper $H$-invariant subspace. It is nondegenerate by Corollary~\ref{cor:red}.

Write the signature of $W$ as $(a,b)$. Since $W$ contains isotropic lines, $a,b\ge1$. Viewed inside $\partial_\infty\HH(W)$, the set $\Lambda$ is still compact and negative, since the form and all pairings are obtained by restricting the ambient form to $W$. Since $\Lambda\cong\mathbb S^{p-1}$, Lemma~\ref{lem:acausal-dim} gives $p-1=\dim\Lambda\le a-1$. On the other hand, $a\le p$, since $W\subset\R^{p,q}$. Hence $a=p$, so $W^\perp$ is negative definite. Since $W$ is proper, $1\le b\le q-1$.

Every orthogonal lift of an element of $\rho(\Gamma)$ is block diagonal with respect to $V=W\oplus W^\perp$. Its $W^\perp$-block preserves a negative definite form, so all eigenvalues of this block have modulus one. Since $W$ has signature $(p,b)$ with $p>b$, the $W$-block also has eigenvalues of modulus one. Consequently the two largest eigenvalue moduli, and hence the first Jordan gap, are unchanged upon restriction to $W$. The map $\xi$, now viewed as taking values in $\Pj(W)$, remains transverse and dynamics-preserving. The Jordan-gap characterization therefore shows that $\rho_W$ is $P_1$-Anosov, with the same boundary map. For $b\ge2$, its boundary map is negative, so Theorem~\ref{thm:dgk-characterization} shows that $\rho_W$ is $\HH^{p,b-1}$-convex cocompact. For $b=1$, the same conclusion is the usual rank-one characterization of convex cocompact representations by the Anosov property. Its kernel is finite and its image is discrete. Since $W$ is nondegenerate, $\partial_\infty\HH(W)$ is an embedded submanifold of the ambient boundary. Its conformal structure is induced by the restriction of the same bilinear form to $W$. Therefore $\Lambda\subset\Pj(W)$ is a $C^1$ submanifold of $\partial_\infty\HH(W)$, and the tangent form on $T\Lambda$ is the same as the ambient tangent form. In particular, its radical rank is unchanged.

Suppose that $\Zcl(\rho_W(\Gamma))$ were reducible on $W$, and let $U<W$ be a nonzero proper invariant subspace. Then $U$ is invariant under $\rho(\Gamma)$ and hence under its ambient Zariski closure $H$. Corollary~\ref{cor:red} therefore shows that $U$ is nondegenerate, so $W$ admits a nontrivial invariant orthogonal decomposition. Lemma~\ref{lem:directsum} would then place $\Lambda$ in one proper summand, contradicting the definition of $W$. Thus the projected Zariski closure acts irreducibly.

\end{proof}

Consequently, the reducible case in signature $(p,q)$ reduces to the same rigidity problem in a strictly smaller negative index. If $b=1$, or if $b=2$ by the theorem of Glorieux--Monclair \cite[Theorem~1.3]{GM24}, then $\rho(\Gamma)$ preserves a totally geodesic copy of $\HH^p$.

\section{Excluding the Weyl wall and proof of Theorem~\ref{thm:intro-smallwall}}\label{sec:smallwall}

\subsection{Restricted weights and conventions}

Let $H<\SO(V)$ be a connected reductive group. Following the definition of the limit cone for an arbitrary subsemigroup in Subsection~\ref{sec:projections}, we write $\limcone_H^{\SO(V)}$ for the limit cone of $H$ computed in the ambient group $\SO(V)$; see also \cite[Definition~2.10]{BK25}. This is generally different from the intrinsic cone $\limcone_H^H=\mathfrak a_H^+$, and it should also be distinguished from the cone $\limcone_\Delta^{\SO(V)}$ of a Zariski-dense discrete subgroup $\Delta<H$.

Let $A_H<H$ be a maximal real split torus and put $\mathfrak a=\operatorname{Lie}(A_H)$. If $V$ is a finite-dimensional real algebraic $H$-module, write
$$
V_\chi =\{v\in V:X\cdot v=\chi(X)v\text{ for every }X\in\mathfrak a\},
\qquad
\Omega(V)=\{\chi\in\mathfrak a^*:V_\chi\ne0\}.
$$

Since $A_H$ is split, its action is diagonalizable over $\mathbb R$, and hence
$$
V=\bigoplus_{\chi\in\Omega(V)}V_\chi.
$$

The elements of $\Omega(V)$ are the \emph{restricted weights} of $V$. In particular,
$$
V_0:=\{v\in V:X\cdot v=0\text{ for every }X\in\mathfrak a\}=V^{\mathfrak a}.
$$

For the adjoint representation of $H$ on $\mathfrak h$, write
$$
\mathfrak h=\mathfrak h_0\oplus
\bigoplus_{\beta\in\Sigma}\mathfrak h_\beta,
\qquad
\mathfrak h_\beta
=\{Y\in\mathfrak h:[X,Y]=\beta(X)Y
\text{ for every }X\in\mathfrak a\},
$$
where $\Sigma:=\{\beta\in\mathfrak a^*\setminus\{0\}:\mathfrak h_\beta\ne0\}$ is the restricted root system and $m_\beta=\dim_{\mathbb R}\mathfrak h_\beta$ is the multiplicity of $\beta$. The restricted root system may be nonreduced. Moreover, $\mathfrak h_0=\mathfrak z_{\mathfrak h}(\mathfrak a)$. If $\mathfrak h=\mathfrak k\oplus\mathfrak p$ is a compatible Cartan decomposition, $\mathfrak a\subset\mathfrak p$, and $\mathfrak m=\mathfrak z_{\mathfrak k}(\mathfrak a)$, then $\mathfrak h_0=\mathfrak m\oplus\mathfrak a$. Thus $\mathfrak a\subset\mathfrak h_0$, but equality need not hold even in the adjoint representation. For an arbitrary module, the root-space action satisfies $\mathfrak h_\beta V_\chi\subset V_{\chi+\beta}$, where the right-hand side is understood to be zero if $\chi+\beta$ is not a weight. The restricted Weyl group $W(\mathfrak h,\mathfrak a)=N_H(\mathfrak a)/Z_H(\mathfrak a)$ acts on $\mathfrak a^*$ and preserves both $\Sigma$ and $\Omega(V)$, together with their multiplicities.

For example, in the standard representation of $\SO_{p,q}^{\circ}$, with $p>q$, one may choose coordinates so that
$$
\Omega(V)=\{0,\pm\varepsilon_1,\ldots,\pm\varepsilon_q\},
\qquad
\dim V_{\pm\varepsilon_i}=1,
\qquad
\dim V_0=p-q.
$$

For $X\in\mathfrak a$, put
$$
n_V(X)=\sum_{\chi(X)>0}\dim_{\mathbb R}V_\chi.
$$

For a complex module $W$, we use the same notation $n_W(X)$, with complex dimensions in place of real dimensions. We call $X$ \emph{weight-regular} if $\chi(X)\ne0$ for every nonzero $\chi\in\Omega(V)$. If the representation is orthogonal, invariance of the form gives $B(V_\chi,V_\eta)=0$ unless $\chi+\eta=0$, and $B$ induces a perfect pairing between $V_\chi$ and $V_{-\chi}$. Thus the nonzero restricted weights occur in opposite pairs with equal multiplicities. We use the highest-weight theorem \cite[Theorem~5.5]{Knapp} and the restricted-weight formalism of \cite[Lemma~8.48]{Knapp}. For the restricted-root decomposition and the restricted Weyl group, see \cite[Chapter~VI, Sections~1--3]{Helgason}; the classification and root multiplicities are recorded in \cite[Chapter~X, Table~VI]{Helgason}.

We also fix some standard representation-theoretic notation. For a real module $V$, let $V_{\mathbb C}=V\otimes_{\mathbb R}\mathbb C$ denote its complexification; for a complex module $W$, let $W_{\mathbb R}$ denote its underlying real module and $W^\sigma$ its conjugate module. If $H=H_1\times\cdots\times H_s$ (or $\mathfrak h=\mathfrak h_1\oplus\cdots\oplus\mathfrak h_s$) and $W_i$ is a module for the $i$-th factor, then $W_1\boxtimes\cdots\boxtimes W_s$ denotes the \emph{external tensor product}: its underlying vector space is $W_1\otimes\cdots\otimes W_s$, and each factor acts on the corresponding tensor factor. Every finite-dimensional irreducible complex module for a direct sum of complex semisimple Lie algebras is of this form. We write $\operatorname{Sym}^mW$ and $\Lambda^mW$ for the symmetric and exterior powers; if $W$ carries a nondegenerate invariant symmetric bilinear form, $\operatorname{Sym}^m_0W$ denotes the traceless, or harmonic, symmetric power, namely the kernel of the contraction $\operatorname{Sym}^mW\to\operatorname{Sym}^{m-2}W$.

An irreducible real module $V$ is of real, complex, or quaternionic type according as $\operatorname{End}_H(V)$ is isomorphic to $\mathbb R$, $\mathbb C$, or $\mathbb H$; correspondingly, $V_{\mathbb C}$ is irreducible, is $W\oplus W^\sigma$ with $W\not\cong W^\sigma$, or is the direct sum of two isomorphic copies of a self-conjugate complex irreducible. For a self-conjugate irreducible complex module $W$, its Frobenius--Schur sign is $+1$ if it admits an equivariant conjugate-linear map $J$ with $J^2=\id$, and $-1$ if it admits one with $J^2=-\id$. In the latter case, $J$, viewed as a real-linear endomorphism of $W_{\mathbb R}$, is a complex structure. If $I$ denotes multiplication by $i$, then $IJ=-JI$, and $I,J,K:=IJ$ define a quaternionic structure on $W_{\mathbb R}$. For self-conjugate factors, Frobenius--Schur signs multiply under external tensor products. Throughout this section, the \emph{local type} of a group means the isomorphism type of its Lie algebra.

\subsection{Reduction to real rank at most two}\label{sec:smallwall-reduction}

The following lemma bounds the number of simple factors. Its proof distinguishes real, complex, and quaternionic representation types, since the complexification of an irreducible real module can have either one or two irreducible constituents, counted with multiplicity.

\begin{lemma}\label{lem:real-factor}

Let $\mathfrak h$ be a noncompact real semisimple Lie algebra and let $V$ be a faithful irreducible real $\mathfrak h$-module. Let $\mathfrak a$ be the Lie algebra of a maximal real split torus. Assume that
$$
n_V(A)\le2 \quad \text{for all} \quad A\in\mathfrak a.
$$

Then either $\mathfrak h$ is noncompact simple, or
$$
\mathfrak h\cong\mathfrak{sl}_2(\mathbb R)\oplus
\mathfrak{sl}_2(\mathbb R),
\qquad
V\cong\mathbb R^2\boxtimes\mathbb R^2.
$$

In the exceptional case $\dim_{\mathbb R}V=4$, and every invariant nondegenerate symmetric form has signature $(2,2)$. In particular, if $V=\mathbb R^{p,3}$ with $p>3$, then $\mathfrak h$ is noncompact simple.
\end{lemma}

\begin{proof}
Write
$$
\mathfrak h=\mathfrak h_1\oplus\cdots\oplus\mathfrak h_s
$$
as a sum of real simple ideals. If $W$ is an irreducible constituent of $V_{\mathbb C}$, the real, complex, and quaternionic trichotomy recalled above gives
$$
V_{\mathbb C}\cong W
\quad\hbox{or}\quad
V_{\mathbb C}\cong W\oplus W^\sigma .
$$

The second alternative includes both complex type and quaternionic type; in the latter case the two summands are isomorphic. Put $m=1$ or $2$ according to the two alternatives. Conjugation does not alter restricted weights because they are real on $\mathfrak a$, and $(V_\chi)_{\mathbb C}=(V_{\mathbb C})_\chi$. Thus $\dim_{\mathbb R}V_\chi=m\dim_{\mathbb C}W_\chi$.

The classification of irreducible modules for a direct sum gives $W=W_1\boxtimes\cdots\boxtimes W_s$, where $W_i$ is irreducible for $(\mathfrak h_i)_{\mathbb C}$. If $\mathfrak h_i$ is a complex simple Lie algebra viewed as real, then $(\mathfrak h_i)_{\mathbb C}$ is a sum of two simple ideals. Every $W_i$ is nontrivial by faithfulness; put $d_i=\dim_{\mathbb C}W_i\ge2$. If $\mathfrak h_i$ is noncompact and $A_i$ belongs to the Lie algebra of a maximal real split torus of this factor, with all other components zero, then
\begin{equation}\label{eq:real-factor-count}
n_V(A_i)=m\,n_{W_i}(A_i)\prod_{j\ne i}d_j .
\end{equation}

The Lie algebra of the maximal split torus acts nontrivially on $W_i$, since otherwise the ideal generated by it would act trivially. As the action is diagonalizable with real eigenvalues and has trace zero, a suitable $A_i$ has both positive and negative eigenvalues; hence $n_{W_i}(A_i)\ge1$. The wall bound and \eqref{eq:real-factor-count} imply
\begin{equation}\label{eq:factor-budget}
m\prod_{j\ne i}d_j\le2.
\end{equation}

Thus $s\le2$. If $s=2$ and $\mathfrak h_1$ is noncompact, then $m=1$ and $d_2=2$.

Suppose first that $\mathfrak h_2$ is compact. Its nontrivial two-dimensional complex module forces $\mathfrak h_2\cong\mathfrak{su}_2$, and $W_2$ is quaternionic. Since external tensor factors are unique and $W_2\cong W_2^\sigma$, self-conjugacy of $W_1\boxtimes W_2=V_{\mathbb C}$ forces $W_1\cong W_1^\sigma$. Frobenius--Schur signs then show that $W_1$ is quaternionic: the sign of the tensor product is $+1$, while that of $W_2$ is $-1$. Its equivariant antilinear quaternionic structure preserves every restricted-weight space, so each of those spaces has even complex dimension. Hence $n_{W_1}(A_1)\ge2$, and \eqref{eq:real-factor-count} gives $n_V(A_1)=2n_{W_1}(A_1)\ge4$, a contradiction.

Both factors are therefore noncompact. Applying \eqref{eq:factor-budget} in both directions gives $m=1$ and $d_1=d_2=2$. The real structure forces the two factors to be self-conjugate. If the complexification of $\mathfrak h_i$ is simple, its nonzero two-dimensional module embeds it into $\mathfrak{sl}_2(\mathbb C)$; hence the complexification is $\mathfrak{sl}_2(\mathbb C)$, and noncompactness gives $\mathfrak h_i\cong\mathfrak{sl}_2(\mathbb R)$. If $\mathfrak h_i$ is complex simple viewed as real, then its complexification is a sum of two ideals; a two-dimensional constituent is supported on one ideal and conjugation moves it to the other, so it is not self-conjugate. Consequently $V\cong\mathbb R^2\boxtimes\mathbb R^2$. If $\omega_i$ is the invariant area form on the $i$-th factor, then $\omega_1\otimes\omega_2$ is symmetric of signature $(2,2)$. This tensor module is absolutely irreducible, so every invariant symmetric form is a scalar multiple of it. This proves the lemma.
\end{proof}

\begin{proposition}\label{prop:general-wall-rank}
Let $p>q\ge2$, and let $H<\SO_{p,q}$ be connected and reductive, acting irreducibly on $V=\R^{p,q}$. If $\limcone_H^{\SO_{p,q}}\subset F_{q-1}$, then $n_V(A)\le q-1$ for every $A\in\mathfrak a$, and $\rank_{\R}H\le q-1$. If at least one of $p$ and $q$ is odd, then the Lie algebra of $H$ is semisimple.
\end{proposition}

\begin{proof}
For $A\in\mathfrak a$, the equality $\lambda_q(\exp A)=0$ says exactly that at most $q-1$ eigenvalue moduli of $\exp A$ exceed $1$, counted with multiplicity. This is the inequality $n_V(A)\le q-1$.

The nonzero restricted weights span $\mathfrak a^*$: otherwise a nonzero element of $\mathfrak a$ annihilated by all the weights would act trivially on the faithful module $V$. For a weight-regular $A$, orthogonality pairs the nonzero weights as $\{\chi,-\chi\}$, and exactly one member of each pair is positive on $A$. Hence the positive weights span $\mathfrak a^*$. Since every positive weight space is nonzero and their total dimension is at most $q-1$, there are at most $q-1$ positive weights. Therefore $\dim\mathfrak a\le q-1$. The last assertion is Lemma~\ref{lem:semisimple-irreducible}.
\end{proof}

\begin{proposition}[Rank reduction in index three]\label{prop:rankred}
Let $H<\SO_{p,3}$, $p>3$, be connected, noncompact, reductive, and irreducible on $V=\mathbb R^{p,3}$. If $\limcone_H^{\SO_{p,3}}\subset F_2$, then the Lie algebra of $H$ is noncompact simple and has real rank at most two.
\end{proposition}
\begin{proof}
Proposition~\ref{prop:general-wall-rank} gives
\begin{equation}\label{eq:small-wall-count}
n_V(A)\le2\qquad(A\in\mathfrak a).
\end{equation}

It also gives $\rank_{\R}H\le2$ and shows that $\mathfrak h$ is semisimple. The exceptional alternative in Lemma~\ref{lem:real-factor} is four-dimensional, whereas $\dim V=p+3\ge7$; hence $\mathfrak h$ is noncompact simple.
\end{proof}

\begin{remark}[What is special to index three]\label{rem:section-six-general-q}
The restricted-weight formalism and the classification results below are not tied to an ambient signature: they classify irreducible orthogonal modules satisfying the exceptional numerical bound $n_V\le2$. What is special to $\SO_{p,3}$ is that its last wall $F_2$ gives exactly this bound. In $\SO_{p,q}$ the last wall $F_{q-1}$ gives only $n_V\le q-1$, as Proposition~\ref{prop:general-wall-rank} shows. The analogue of Theorem~\ref{thm:intro-smallwall} with $3$ replaced by $q$ is false for every $q\ge4$: the irreducible adjoint representation, equipped with the form $-B_{\mathfrak{so}_{q,1}}$, gives $\Ad(\SO_0(q,1))<\SO_{\binom{q}{2},q}$, and a regular split element has ambient Jordan projection $(\underbrace{t,\ldots,t}_{q-1},0)$, where $t>0$. Thus $\limcone_{\Ad(\SO_0(q,1))}^{\SO_{\binom{q}{2},q}}$ lies in the last wall $F_{q-1}$. This does not contradict the geometric conjecture: the image of this adjoint representation contains no proximal elements, since its highest restricted-weight space has dimension $q-1$. Proximality alone does not restore a general exclusion: the image $\Ad(\SO_0(4,2))<\SO_{8,7}$ contains proximal elements and its cone lies in the last wall of $\SO_{8,7}$. On the other hand, the stronger face $F_2$ in any larger signature still gives the bound $n_V\le2$, so the representation-theoretic statements below continue to apply once their stated simplicity and irreducibility hypotheses are available.
\end{remark}

We now classify the irreducible orthogonal modules satisfying $n_V\le2$. The final step will be to compare their signatures with $(p,3)$.

\begin{theorem}[Low-rank wall classification]\label{thm:classification}
Let $H$ have noncompact simple real Lie algebra of real rank at most two, and let $V$ be a nontrivial irreducible real orthogonal $H$-module. If $n_V(A)\le2$ for a weight-regular $A$ in the Lie algebra of a maximal real split torus, then, up to local isomorphism and replacing the invariant form by its negative, $V$ is one of the following:
$$
\begin{array}{c|c|c}
\text{local type of }H&V&\text{signature}\\
\hline
\SL_2(\mathbb R)&\operatorname{Sym}^2\mathbb R^2,\
\operatorname{Sym}^4\mathbb R^2&(2,1),\ (3,2)\\
\SO_{n,1}^{\circ},\ n\ge3&\mathbb R^{n,1}&(n,1)\\
\SO_{3,1}^{\circ}&\mathfrak{so}_{3,1}&(3,3)\\
\SU_{n,1},\ n\ge2&(\mathbb C^{n,1})_{\mathbb R}&(2n,2)\\
\SO_{n,2}^{\circ},\ n\ge3&\mathbb R^{n,2}&(n,2)
\end{array}
$$

In particular, none has signature $(p,3)$ with $p>3$.
\end{theorem}

The next two subsections prove Theorem~\ref{thm:classification}. Throughout, $A_H$ denotes a maximal real split torus, $\mathfrak a$ its Lie algebra, and $A_0\in\mathfrak a$ a fixed weight-regular element.

\subsection{Real-rank-one groups}\label{sec:rankone-classification}

\begin{lemma}\label{lem:rankone-realification} 
Let $G$ be a connected real semisimple group, let $A_G<G$ be a maximal real split torus with Lie algebra $\mathfrak a$, and let $V$ be an irreducible real $G$-module. Fix a weight-regular $A_0\in\mathfrak a$ and put
$$
N_{A_0}(V):=n_V(A_0)=\sum_{\chi(A_0)>0}\dim_{\mathbb R}V_\chi.
$$

For a complex module $W$, put $N_{\mathbb C,A_0}(W):=n_W(A_0)$, so the dimensions in this count are taken over $\mathbb C$. Then exactly one of the following occurs.
\begin{enumerate}
\item $V_{\mathbb C}$ is an irreducible complex module $W$, and $N_{A_0}(V)=N_{\mathbb C,A_0}(W)$. 
\item $V_{\mathbb C}=W\oplus W^\sigma$, where $W^\sigma$ is the conjugate of $W$ (possibly isomorphic to $W$). In this case $N_{A_0}(V)=2N_{\mathbb C,A_0}(W)$.
\end{enumerate}
The second alternative comprises both complex type and quaternionic type. 
\end{lemma}
\begin{proof}
This is the usual real/complex/quaternionic trichotomy for irreducible real representations. Since every restricted weight is real on $\mathfrak a$, conjugation does not change it. Moreover $(V_\chi)_{\mathbb C}=(V_{\mathbb C})_\chi$, so comparison of dimensions gives the two formulas in the statement.
\end{proof}

In the rank-one calculations below, we choose $A_0$ in the positive half-line of $\mathfrak a$ and write $N(V):=n_V(A_0)$ and $N_{\mathbb C}(W):=n_W(A_0)$.

We shall also use the following standard highest-weight observation. Let $\mathfrak g_{\mathbb C}$ be complex semisimple, let $h$ lie in the closed dominant Weyl chamber, and let $W_\lambda$ be the irreducible module of highest weight $\lambda$. Fix a nonzero highest-weight vector $v_\lambda$ and, for each root $\beta$, choose raising and lowering operators $e_\beta$ and $f_\beta$. Put $\mathfrak l_h=[Z_{\mathfrak g_{\mathbb C}}(h),Z_{\mathfrak g_{\mathbb C}}(h)]$ and $E_{\max}(h,W_\lambda)=\{v\in W_\lambda:h\cdot v=\lambda(h)v\}$. Since $h$ is dominant, $\lambda(h)$ is the largest $h$-eigenvalue on $W_\lambda$. The space $E_{\max}(h,W_\lambda)$ is the irreducible $\mathfrak l_h$-module generated by $v_\lambda$, and its Dynkin labels are the labels of $\lambda$ at the simple roots $\alpha$ satisfying $\alpha(h)=0$. See \cite[Theorem~5.104 and Proposition~5.105]{Knapp}. We shall also repeatedly use the following root-string consequence of the highest-weight theorem \cite[Theorem~5.5]{Knapp}:
$$
\langle\lambda,\beta^\vee\rangle>0 \quad\Longrightarrow\quad
\lambda-\beta\text{ is a weight of }W_\lambda .
$$

\begin{proposition}[Rank-one wall classification]\label{prop:rankone-small-wall-complete} 
Let $H$ be a connected noncompact simple real group of real rank one and let $V$ be a nontrivial irreducible real orthogonal $H$-module. Suppose that
$$
N(V)=\sum_{\chi>0}\dim_{\mathbb R}V_\chi\le2.
$$

Up to local isomorphism, meaning up to the Lie algebra of $H$, the complete list of underlying $H$-modules is as follows; every invariant nondegenerate symmetric form has the indicated unordered signature:
\begin{equation}\label{eq:rankone-survivors}
\begin{array}{c|c|c|c}
\text{local type of }H&V&N(V)&\text{signature}\cr
\hline
\SL_2(\mathbb R)&\operatorname{Sym}^2\mathbb R^2&1&(2,1)\cr
\SL_2(\mathbb R)&\operatorname{Sym}^4\mathbb R^2&2&(3,2)\cr
\SO_{n,1}^{\circ},\ n\ge3&\mathbb R^{n,1}&1&(n,1)\cr
\SO_{3,1}^{\circ}&\mathfrak{so}_{3,1}&2&(3,3)\cr
\SU_{n,1},\ n\ge2&(\mathbb C^{n,1})_{\mathbb R}&2&(2n,2)\cr
\Sp_{1,1}\cong\Spin_{4,1}&\mathbb R^{4,1}&1&(4,1).
\end{array}
\end{equation}

The last row is the $n=4$ instance of the orthogonal row; we repeat it to record the $\Sp_{1,1}$ case. In particular, no such module has signature $(p,3)$ with $p>3$.
\end{proposition}

\begin{proof}
The list of noncompact simple real Lie algebras of real rank one is
$$
\mathfrak{sl}_2(\mathbb R),\quad
\mathfrak{so}_{n,1},\quad
\mathfrak{su}_{n,1},\quad
\mathfrak{sp}_{n,1},\quad
\mathfrak f_{4(-20)};
$$
we take $n\ge3$ in the orthogonal family, $n\ge2$ in the unitary family, and $n\ge1$ in the quaternionic family; the overlap arising from the exceptional isomorphism $\mathfrak{sp}_{1,1}\cong\mathfrak{so}_{4,1}$ is retained to make the quaternionic family explicit. The restricted-root data used below agree with the rank-one rows of \cite[Ch.~X, Table~VI]{Helgason}. We give the weight calculations, including all low-rank exceptions.

\smallskip \noindent\emph{1. The split group $\SL_2(\mathbb R)$.} Every complex irreducible is $\operatorname{Sym}^m\mathbb C^2$, and the corresponding irreducible real module is $\operatorname{Sym}^m\mathbb R^2$. Its weights are $m,m-2,\ldots,-m$. The invariant bilinear form is symmetric exactly when $m=2k$ is even, and then $N(V)=k$. Thus $N(V)\le2$ gives $k=1,2$. In the usual monomial basis the form pairs the $j$-th and $(2k-j)$-th basis vectors, with alternating signs. The $k$ pairs not involving the middle basis vector span hyperbolic planes, while the restriction of the form to the middle line is nondegenerate. Hence the unordered signature is $(k+1,k)$, proving the first two rows of \eqref{eq:rankone-survivors}.

\smallskip \noindent\emph{2. The groups $\SO_{n,1}^{\circ}$.} First assume $n\ge4$. Put $N=n+1$, so that $\mathfrak g_{\mathbb C}=\mathfrak{so}_N(\mathbb C)$, and choose coordinates in which $e_1(h)=1$ and $e_j(h)=0$ for $j\ge2$. The vector representation has $h$-weights $1,0^{(n-1)},-1$. If $\lambda=(\lambda_1,\lambda_2,\ldots)$ is a dominant highest weight, $E_{\max}(h,W_\lambda)$ is the irreducible $\mathfrak{so}_{N-2}(\mathbb C)$-module of highest weight $(\lambda_2,\lambda_3,\ldots)$.

For $n\ge6$, every nontrivial complex representation of $\mathfrak{so}_{n-1}(\mathbb C)$ has dimension at least $4$. Consequently a module with $N_{\mathbb C}\le2$ must have $\lambda=j\omega_1$. It is the traceless harmonic power $\operatorname{Sym}^j_0\mathbb C^N$, namely the kernel of contraction on $\operatorname{Sym}^j\mathbb C^N$. For $j\ge2$, its two largest $h$-eigenvalues have multiplicities $1$ and $n-1$: if $e_+$ spans the $h$-weight-one line and $u$ ranges over the $(n-1)$-dimensional zero-weight space, the corresponding eigenspaces are spanned by $e_+^j$ and by the vectors $e_+^{j-1}u$, respectively; all these vectors lie in the kernel of contraction. Thus $N_{\mathbb C}\ge n$. Only $j=1$, the vector module, remains.

The cases $n=4$ and $n=5$, in which the centralizer has low rank, require separate consideration. For $n=4$ the root system is $B_2$, and we write $\lambda=a\omega_1+b\omega_2$. The module $E_{\max}(h,W_\lambda)$, regarded as an $\mathfrak{so}_3(\mathbb C)$-module has dimension $b+1$. Hence $b\ge2$ is impossible. If $b=1$ and $a\ge1$, $f_{e_1-e_2}v_\lambda$ is nonzero, lies outside $E_{\max}(h,W_\lambda)$, and has positive $h$-eigenvalue. Thus the only possibility with $b\ne0$ is $\omega_2$, the four-dimensional complex spin module, for which $N_{\mathbb C}=2$; when $b=0$, the harmonic-power calculation above excludes $a\ge2$. For $n=5$, write a dominant $D_3$ weight as $\lambda_1\ge\lambda_2\ge|\lambda_3|$. The $\mathfrak{so}_4(\mathbb C)$-module $E_{\max}(h,W_\lambda)$ has dimension $(\lambda_2+\lambda_3+1)(\lambda_2-\lambda_3+1)$. It has dimension at most two only when it is trivial or when $(\lambda_2,\lambda_3)=(\tfrac12,\pm\tfrac12)$. In the latter case, $\lambda_1\in\tfrac12+\mathbb Z$, so $\lambda_1>\tfrac12$ implies $\lambda_1\ge\tfrac32$; lowering by $e_1-e_2$ then gives a nonzero vector outside $E_{\max}(h,W_\lambda)$ with positive $h$-eigenvalue $\lambda_1-1$. The only new possibilities are therefore the two complex half-spin modules, each with $N_{\mathbb C}=2$. If $E_{\max}(h,W_\lambda)$ is trivial, then $\lambda=j\omega_1$, and the same harmonic-power calculation leaves only $j=1$.

Neither spin exception gives a real module with $N(V)\le2$. Indeed
$$
\Spin_{4,1}\cong\Sp_{1,1},\qquad
\Spin_{5,1}\cong \SU^*(4)\cong\SL_2(\mathbb H),
$$
and in each case the indicated complex spin module is the defining quaternionic module (or its dual) under the corresponding accidental isomorphism. It has an equivariant antilinear operator $J$ with $J^2=-\id$. Hence its irreducible underlying real module has complexification equal to two copies of the spin module, and Lemma~\ref{lem:rankone-realification} gives $N(V)=4$. We have proved that, for $n\ge4$, the vector module is the only irreducible real module satisfying $N(V)\le2$.

It remains to treat the case $n=3$. Here
$$
\mathfrak{so}_{3,1}\cong\mathfrak{sl}_2(\mathbb C)_{\mathbb R},\qquad
\mathfrak g_{\mathbb C}=\mathfrak{sl}_2(\mathbb C)\oplus
\mathfrak{sl}_2(\mathbb C).
$$

Write $W_{a,b}=\operatorname{Sym}^a\mathbb C^2\boxtimes\operatorname{Sym}^b\mathbb C^2$. Let $H_0=\operatorname{diag}(1,-1)$. Real conjugation interchanges the two factors, and $h=(H_0,H_0)$. Let $P(a,b)$ be the number, with multiplicity, of positive sums of a weight of $\operatorname{Sym}^a$ and a weight of $\operatorname{Sym}^b$. If $a=b$, the module has a real form and $N(V)=P(a,a)$. Now $P(1,1)=1$, whereas $P(a,a)\ge3$ for $a\ge2$: the highest-weight vector and the two vectors obtained by applying the simple-root lowering operator in the first and second factors, respectively, all have positive $h$-eigenvalue. If $a\ne b$, an irreducible real module has complexification $W_{a,b}\oplus W_{b,a}$ and hence $N(V)=2P(a,b)$. Thus $N(V)\le2$ forces $P(a,b)=1$. Assuming $a>b$, this happens only for $(a,b)=(1,0)$ or $(2,0)$: if $b=0$, then $P(a,0)=\lceil a/2\rceil$, while if $a>b\ge1$, the highest-weight vector and the vector obtained by applying the simple-root lowering operator in the first factor both have positive $h$-eigenvalue. The first possibility is symplectic, not orthogonal. More generally, $W_{a,b}$ has a symmetric invariant complex bilinear form exactly when $a+b$ is even. For $a\ne b$, the modules $W_{a,b}$ and $W_{b,a}$ are nonisomorphic and each is self-dual, so an invariant bilinear form on their direct sum has no cross terms. Thus the only orthogonal case with $a\ne b$ is $(2,0)$, up to interchanging the factors. The surviving real modules are therefore the vector representation of signature $(3,1)$, whose complexification is $W_{1,1}$, and the adjoint representation on $\mathfrak{so}_{3,1}$, of signature $(3,3)$, whose complexification is $W_{2,0}\oplus W_{0,2}$.

\smallskip \noindent\emph{3. The groups $\SU_{n,1}$, $n\ge2$.} Put $N=n+1$ and take $\mathfrak g_{\mathbb C}=\mathfrak{sl}_N(\mathbb C)$ and $h=\operatorname{diag}(1,0,\ldots,0,-1)$. For a highest weight $\lambda=\sum_{i=1}^{N-1}a_i\omega_i$, the largest $h$-eigenvalue is $S=\lambda(h)=\sum_i a_i$. The corresponding eigenspace is irreducible for the semisimple part $\mathfrak l_h\cong\mathfrak{sl}_{N-2}(\mathbb C)$ of $Z_{\mathfrak g_{\mathbb C}}(h)$ and has Dynkin labels $(a_2,\ldots,a_{N-2})$.

If $N\ge5$, a nontrivial such top module has dimension at least $N-2\ge3$, so all middle labels vanish. If both endpoint labels $a_1,a_{N-1}$ are nonzero, the top vector and its lowerings by the two endpoint simple roots give three linearly independent vectors with positive $h$-eigenvalues. If only one endpoint label is nonzero, the module is $\operatorname{Sym}^a\mathbb C^N$ or its dual. For $a=1$ this is the defining module and $N_{\mathbb C}=1$; for $a=2$ the positive eigenspaces have dimensions $1$ and $N-2$, and for $a\ge3$ the first three vectors in the endpoint $\mathfrak{sl}_2$-string have positive $h$-eigenvalues. Thus only the defining module and its dual have $N_{\mathbb C}\le2$.

For $N=4$, one has $\dim E_{\max}(h,W_\lambda)=a_2+1$. When $a_2=0$, the endpoint argument from the preceding paragraph, which applies also for $N=4$, leaves only the defining module and its dual. If $a_2=1$ and at least one endpoint label is nonzero, applying the corresponding endpoint lowering operator to $v_\lambda$ gives a nonzero vector of positive $h$-eigenvalue outside $E_{\max}(h,W_\lambda)$; if $a_2\ge2$, then $\dim E_{\max}(h,W_\lambda)\ge3$. Thus the only additional highest weight is $\omega_2$, and $\Lambda^2\mathbb C^4$ has exactly two weights $\mu$ satisfying $\mu(h)>0$. For $N=3$, if both Dynkin labels are nonzero, the highest-weight vector and its images under the two simple-root lowering operators are three linearly independent $h$-eigenvectors with positive eigenvalues. If $\lambda=a\omega_1$ or $a\omega_2$, the corresponding simple-root string shows that only $a=1,2$ can occur. Consequently,
$$
\begin{array}{c|c}
N&\text{complex modules with }N_{\mathbb C}\le2\cr
\hline
3&\mathbb C^3,(\mathbb C^3)^*,\operatorname{Sym}^2\mathbb C^3,
\operatorname{Sym}^2(\mathbb C^3)^*\cr
4&\mathbb C^4,(\mathbb C^4)^*,\Lambda^2\mathbb C^4\cr
\ge5&\mathbb C^N,(\mathbb C^N)^*.
\end{array}
$$

For the real form $\SU_{n,1}$, conjugation sends a complex module to its dual. Indeed, the Hermitian form identifies $\overline{\mathbb C^N}$ with $(\mathbb C^N)^*$ and hence identifies the conjugate of every Schur module with its dual. The defining module and its dual therefore occur together in the complexification of the irreducible real module $(\mathbb C^{n,1})_{\mathbb R}$; it has $N(V)=2$ and the real part of the Hermitian form has signature $(2n,2)$. The two symmetric-square modules for $\SU_{2,1}$ likewise combine, and hence have $N(V)=4$. Finally, equip $\mathbb C^4$ with its Hermitian form of signature $(3,1)$ and choose a pseudo-orthonormal basis $e_1,e_2,e_3,e_4$, with $e_4$ negative and volume $e_{1234}$. On $\Lambda^2\mathbb C^4$ the equivariant conjugate-linear Hodge operator satisfies
$$
*e_{12}=e_{34},\qquad *e_{34}=-e_{12},\qquad *^2=-\id.
$$

Every equivariant conjugate-linear automorphism is, by Schur's lemma, of the form $c*$ with $c\in\mathbb C^*$. Its square is $-|c|^2\id$, so none is a real structure; the module is of quaternionic type. The complexification of its irreducible underlying real module is $\Lambda^2\mathbb C^4\oplus\Lambda^2\mathbb C^4$, and hence $N(V)=4$. This proves the unitary row of \eqref{eq:rankone-survivors}.

\smallskip \noindent\emph{4. The groups $\Sp_{n,1}$.} Put $\ell=n+1$. For $\mathfrak g_{\mathbb C}=\mathfrak{sp}_{2\ell}(\mathbb C)$ choose orthogonal coordinates $e_1,\ldots,e_\ell$ such that $e_1(h)=e_2(h)=1$ and $e_j(h)=0$ for $j\ge3$. This gives restricted-root multiplicities $m_\alpha=4(n-1)$ and $m_{2\alpha}=3$. Write a dominant highest weight as $\lambda_1\ge\lambda_2\ge\cdots\ge\lambda_\ell\ge0$. The space $E_{\max}(h,W_\lambda)$ is an irreducible $\mathfrak l_h\cong\mathfrak{sl}_2(\mathbb C)\oplus\mathfrak{sp}_{2\ell-4}(\mathbb C)$-module, with highest weight $\lambda_1-\lambda_2$ on the first factor and $(\lambda_3,\ldots,\lambda_\ell)$ on the second.

Suppose first that $\ell\ge4$. Since the smallest nontrivial module of the second factor has dimension $2(\ell-2)\ge4$, the inequality $\dim E_{\max}(h,W_\lambda)\le2$ forces $\lambda_3=\cdots=\lambda_\ell=0$ and $\lambda_1-\lambda_2\in\{0,1\}$. If $\lambda=(b,b,0,\ldots)=b\omega_2$, then $E_{\max}(h,W_\lambda)$ is one-dimensional. For $b\ge1$, the vector $f_{\alpha_2}v_\lambda$, where $\alpha_2=e_2-e_3$, is a highest-weight vector of highest weight $1$ for the $\mathfrak{sl}_2$-factor of $\mathfrak l_h$ associated with $\alpha_1=e_1-e_2$. The two-dimensional module that it generates lies in the $h$-eigenspace of eigenvalue $2b-1>0$, and hence $N_{\mathbb C}\ge3$. If $\lambda=(b+1,b,0,\ldots)=\omega_1+b\omega_2$, then $\dim E_{\max}(h,W_\lambda)=2$. For $b>0$, the vector $f_{\alpha_2}v_\lambda$ is nonzero, lies outside this eigenspace, and has $h$-eigenvalue $2b>0$. Only $b=0$, namely the defining module $\omega_1$, remains, and it has $N_{\mathbb C}=2$.

For $\ell=3$, the factor $\mathfrak{sp}_2(\mathbb C)\cong\mathfrak{sl}_2(\mathbb C)$ yields one additional family for which $\dim E_{\max}(h,W_\lambda)=2$, namely $\lambda=(b,b,1)$ with $b\ge1$. If $b>1$, then $f_{e_2-e_3}v_\lambda$ is nonzero, lies outside $E_{\max}(h,W_\lambda)$, and has $h$-eigenvalue $2b-1>0$. If $b=1$, so that $\lambda=\omega_3$, then $f_{e_2+e_3}v_\lambda$ is nonzero and has weight $e_1$, hence $h$-eigenvalue $1$. Thus every module in this additional family satisfies $N_{\mathbb C}\ge3$.

For $\ell=2$, the semisimple part $\mathfrak l_h$ of $Z_{\mathfrak g_{\mathbb C}}(h)$ is the $\mathfrak{sl}_2$ associated with $e_1-e_2$. The dominant highest weights for which $\dim E_{\max}(h,W_\lambda)\le2$ are $b\omega_2$ and $\omega_1+b\omega_2$. For $b\omega_2$, lowering by $2e_2$ produces, under $\mathfrak l_h$, a three-dimensional $\mathfrak l_h$-module contained in the $h$-eigenspace of eigenvalue $2b-2$. Indeed, for the lowered weight $\mu=b(e_1+e_2)-2e_2$ one has $\langle\mu,(e_1-e_2)^\vee\rangle=2$, and $e_{e_1-e_2}f_{2e_2}v_\lambda=0$; it therefore generates a three-dimensional $\mathfrak l_h$-module. This $h$-eigenvalue is positive for $b\ge2$, whereas $\omega_2$ itself has $h$-eigenvalues $2,0,-2$, with multiplicities $1,3,1$ and $N_{\mathbb C}=1$. For $\omega_1+b\omega_2$, $E_{\max}(h,W_\lambda)$ has dimension two; if $b>0$, $f_{2e_2}v_\lambda$ is nonzero, lies outside $E_{\max}(h,W_\lambda)$, and has $h$-eigenvalue $2b-1>0$, while $\omega_1$ has $N_{\mathbb C}=2$. We have therefore proved
$$
N_{\mathbb C}(W)\le2
\quad\Longrightarrow\quad
W\cong W_{\omega_1},\quad\text{or}\quad \ell=2\text{ and }W\cong W_{\omega_2}.
$$

The defining complex module $W_{\omega_1}$ has an equivariant quaternionic structure; its underlying irreducible real module has $N(V)=4$ and carries the real part of the quaternionic Hermitian form, of signature $(4n,4)$. When $\ell=2$, the module $W_{\omega_2}$ has a real form: it is the five-dimensional vector module under $\Sp_{1,1}\cong\Spin_{4,1}$, and its form has signature $(4,1)$. This gives precisely the final row of \eqref{eq:rankone-survivors}.

\smallskip \noindent\emph{5. The group $F_{4(-20)}$.} Use the Bourbaki realization
$$
\alpha_1=e_2-e_3,\quad \alpha_2=e_3-e_4,\quad \alpha_3=e_4,
\quad \alpha_4=\tfrac12(e_1-e_2-e_3-e_4).
$$

The dominant cocharacter defining the split torus is
$$
h=(2,0,0,0),\qquad
\alpha_1(h)=\alpha_2(h)=\alpha_3(h)=0,\quad\alpha_4(h)=1.
$$

The eight roots $\beta$ satisfying $\beta(h)=1$ are the half-roots with positive $e_1$-coefficient, while the seven roots satisfying $\beta(h)=2$ are $e_1$ and $e_1\pm e_j$, $2\le j\le4$. Thus the two restricted-root multiplicities are $8$ and $7$, and $\mathfrak l_h$ is of type $B_3$.

Let $\lambda=\sum_{i=1}^4a_i\omega_i$. If one of $a_1,a_2,a_3$ is nonzero, the eigenspace corresponding to the largest $h$-eigenvalue is a nontrivial $B_3$-module and therefore has dimension at least $7$. Otherwise $\lambda=k\omega_4$, $k\ge1$. In these coordinates $\omega_4=e_1$, so the largest $h$-eigenvalue is $2k$. The vector $f_{\alpha_4}v_\lambda$ is nonzero, has $h$-eigenvalue $2k-1$, and is a highest-weight vector of the eight-dimensional spin representation of the semisimple part of $Z_{\mathfrak g_{\mathbb C}}(h)$, which is of type $B_3$, because $-\langle\alpha_4,\alpha_i^\vee\rangle=(0,0,1)$ for $i=1,2,3$. Consequently $N_{\mathbb C}(W_\lambda)\ge1+8=9$. Thus $F_{4(-20)}$ has no nontrivial real module satisfying $N(V)\le2$.

This proves the classification. After fixing an invariant nondegenerate symmetric form $B$, every other such form is of the form $B(T\cdot,\cdot)$ for an invertible $B$-self-adjoint element $T\in\operatorname{End}_H(V)$. Schur's lemma and the representation types determined above show that the space of $B$-self-adjoint endomorphisms is one-dimensional for every module in \eqref{eq:rankone-survivors} except the real adjoint module of $\mathfrak{sl}_2(\mathbb C)$. In this exceptional case the forms are $\operatorname{Re}(cB_{\mathbb C})$, $c\in\mathbb C^*$, where $B_{\mathbb C}$ is the complex Killing form; every one of them has signature $(3,3)$. Thus the signatures in \eqref{eq:rankone-survivors} are fixed up to interchange of the two indices, and none is $(p,3)$ with $p>3$.
\end{proof}

\subsection{Real-rank-two groups}\label{sec:ranktwo-classification}

For groups of real rank two, we use the fact that a faithful orthogonal representation maps each restricted-root space injectively into the corresponding weight space of $\mathfrak{so}(V,B)$.

\begin{proposition}[Rank-two wall classification]\label{prop:ranktwo-smallwall} 
Let $H$ be connected with noncompact simple Lie algebra $\mathfrak h$ of real rank two. Let $\rho:\mathfrak h\longrightarrow\mathfrak{so}(V,B)$ be a faithful irreducible real orthogonal representation. Let $\mathfrak a$ be the Lie algebra of a maximal real split torus, put $\Omega=\Omega(V)$, and write $V=\bigoplus_{\chi\in\Omega}V_\chi$. Suppose that, for a weight-regular $A\in\mathfrak a$,

\begin{equation}\label{eq:ranktwo-budget}
\sum_{\chi(A)>0}\dim V_\chi\le2.
\end{equation}

Then $\mathfrak h\cong\mathfrak{so}_{n,2}$ for some $n\ge3$, and $V$ is its standard $(n+2)$-dimensional representation, up to an automorphism of $\mathfrak h$. In particular, the smaller of the two inertia indices of the invariant form is two.
\end{proposition}

\begin{proof} 
Invariance and nondegeneracy of $B$ give $B(V_\chi,V_\eta)=0$ unless $\eta=-\chi$, and $V_{-\chi}\cong V_\chi^*$. Thus the nonzero restricted weights occur in opposite pairs with equal multiplicities. For weight-regular $A$, \eqref{eq:ranktwo-budget} is equivalent to
\begin{equation}\label{eq:ranktwo-four}
\sum_{\chi\ne0}\dim V_\chi\le4.
\end{equation}
The restricted Weyl group preserves $\Omega$ and all multiplicities.

The restricted root system $\Sigma$ is irreducible of rank two. A nonzero Weyl orbit has at least $3,4,6$ elements in types $A_2,B_2(=C_2),G_2$, respectively. In type $A_2$, an orbit of size three is not stable under multiplication by $-1$, so its negative is a disjoint second orbit. Hence \eqref{eq:ranktwo-four} excludes $A_2$ and $G_2$. In types $B_2,C_2,BC_2$, it implies that $\Omega\setminus\{0\}=S$ is a single four-element Weyl orbit and that every $V_\chi$, $\chi\in S$, is one-dimensional.

For every $\beta\in\Sigma$, faithfulness gives an injective map $\mathfrak h_\beta\hookrightarrow\mathfrak{so}(V,B)_\beta$. In particular, writing $\Omega-\Omega:=\{\chi-\eta:\chi,\eta\in\Omega\}$, one has
\begin{equation}\label{eq:root-difference}
\beta\in\Omega-\Omega.
\end{equation}

Indeed, if $0\ne X\in\mathfrak h_\beta$, then $\rho(X)\ne0$, so $\rho(X)V_\mu\ne0$ for some $\mu\in\Omega$. Since $\rho(X)V_\mu\subset V_{\mu+\beta}$, both $\mu$ and $\mu+\beta$ belong to $\Omega$, and therefore $\beta=(\mu+\beta)-\mu\in\Omega-\Omega$.

Choose coordinates $e_1,e_2$ on $\mathfrak a^*$. Since $\Omega\subset S\cup\{0\}$, condition \eqref{eq:root-difference} implies $\Sigma\subset(S\cup\{0\})-(S\cup\{0\})$. The four-element Weyl orbits satisfying this necessary inclusion are exactly those in the following table:
$$
\begin{array}{c|c|c}
\Sigma&\text{normalization of }\Sigma&\text{candidate four-element orbit }S\\
\hline
B_2&
\{\pm e_i,\ \pm e_1\pm e_2\}&
\{\pm e_1,\pm e_2\},\
\frac12\{\pm e_1\pm e_2\}\\[2mm]
C_2&
\{\pm2e_i,\ \pm e_1\pm e_2\}&
\{\pm e_1,\pm e_2\},\
\{\pm e_1\pm e_2\}\\[2mm]
BC_2&
\{\pm e_i,\ \pm2e_i,\ \pm e_1\pm e_2\}&
\{\pm e_1,\pm e_2\}.
\end{array}
$$

The signs in $\{\pm e_1\pm e_2\}$ vary independently. For $S=c\{\pm e_1,\pm e_2\}$, comparison of $(S\cup\{0\})-(S\cup\{0\})$ with the root systems in the table gives $c=1$ in types $B_2$, $C_2$, and $BC_2$. For $S=c\{\pm e_1\pm e_2\}$, it gives $c=\frac12$ in type $B_2$, $c=1$ in type $C_2$, and no candidate in type $BC_2$.

Use the equivariant identification $\mathfrak{so}(V,B)\cong\Lambda^2V$. For the orbit $S=\frac12\{\pm e_1\pm e_2\}$ in type $B_2$, take $\beta=e_1+e_2$. For the orbit $S=\{\pm e_1,\pm e_2\}$ in type $C_2$ or $BC_2$, take $\beta=2e_1$. In each case $\beta$ is a restricted root, but the only possible contribution to $(\Lambda^2V)_\beta$ is $\Lambda^2V_{\beta/2}$: no two distinct weights in $S\cup\{0\}$ sum to $\beta$. Since $\dim V_{\beta/2}=1$, this exterior square is zero. This contradicts the injectivity of $\mathfrak h_\beta\to\mathfrak{so}(V,B)_\beta$.

We are left with $S=\{\pm e_1,\pm e_2\}$ in type $B_2$ and $S=\{\pm e_1\pm e_2\}$ in type $C_2$. In the first case take $\beta=e_1-e_2$, and in the second take $\beta=2e_1$. The only unordered pairs of distinct weights summing to $\beta$ are, respectively, $\{e_1,-e_2\}$ and $\{e_1+e_2,e_1-e_2\}$. Moreover, $\beta\notin S$, so the zero-weight space contributes nothing. Thus $(\Lambda^2V)_\beta$ is one-dimensional and faithfulness gives
\begin{equation}\label{eq:ranktwo-rootbound}
\dim\mathfrak h_\beta\le1.
\end{equation}

By the classification of real simple Lie algebras and their restricted-root multiplicities \cite[Chapter~X, Table~VI]{Helgason}, the real-rank-two algebras with reduced restricted root system of type $B_2$ or $C_2$ are $\mathfrak{so}_{n,2}$, $n\ge3$, together with $\mathfrak{sp}_{2,2}$ and $\mathfrak{sp}_4(\mathbb C)_{\mathbb R}$. For the last two algebras the restricted root system is $C_2$, and the root $2e_1$ has multiplicity $3$ and $2$, respectively. These multiplicities contradict \eqref{eq:ranktwo-rootbound}. Therefore $\mathfrak h\cong\mathfrak{so}_{n,2}$ for some $n\ge3$. This family includes the accidental isomorphisms
$$
\mathfrak{sp}_4(\mathbb R)\cong\mathfrak{so}_{3,2},\qquad
\mathfrak{su}_{2,2}\cong\mathfrak{so}_{4,2},\qquad
\mathfrak{so}^*(8)\cong\mathfrak{so}_{6,2}.
$$

We now use the $B_2$ normalization for $\mathfrak h\cong\mathfrak{so}_{n,2}$. In the $C_2$ case above, set $f_1=e_1+e_2$ and $f_2=e_1-e_2$, and then rename $f_1,f_2$ as $e_1,e_2$. In these coordinates, the nonzero restricted weights are $\Omega\setminus\{0\}=\{\pm e_1,\pm e_2\}$.

It remains to identify the representation. We first prove that $V_{\mathbb C}$ is irreducible. Suppose instead that $V_{\mathbb C}=W\oplus W^\sigma$. Lemma~\ref{lem:rankone-realification} gives $n_V(A)=2N_{\mathbb C,A}(W)$, where $N_{\mathbb C,A}(W)$ is the total multiplicity of the positive $A$-weights in $W$. Since $W$ is nontrivial, its restricted weights contain a nonzero restricted Weyl orbit. In type $B_2$, every such orbit has cardinality at least four, and the restricted Weyl group contains $-1$. Hence a weight-regular $A$ is positive on at least two weights in this orbit, so $N_{\mathbb C,A}(W)\ge2$ and $n_V(A)\ge4$, contradicting \eqref{eq:ranktwo-budget}.

Choose a Cartan subalgebra $\mathfrak c=\mathfrak a\oplus\mathfrak t$ stable under the Cartan involution, and choose a positive root system compatible with the positive restricted root system. Put $\ell=\lfloor(n+2)/2\rfloor$. In the usual orthogonal coordinates for $\mathfrak{so}_{n+2}(\mathbb C)$, arrange that $\varepsilon_1|_{\mathfrak a}=e_1$, $\varepsilon_2|_{\mathfrak a}=e_2$, and $\varepsilon_j|_{\mathfrak a}=0$ for $j\ge3$. If $\lambda=\lambda_1\varepsilon_1+\cdots+\lambda_\ell\varepsilon_\ell$ is the highest weight of $V_{\mathbb C}$, then $\lambda|_{\mathfrak a}$ is the dominant member of $\Omega\setminus\{0\}=\{\pm e_1,\pm e_2\}$, and hence $\lambda|_{\mathfrak a}=e_1$. Thus $\lambda_1=1$ and $\lambda_2=0$. The dominance inequalities $\lambda_1\ge\lambda_2\ge\cdots\ge\lambda_\ell\ge0$ in type $B_\ell$, or $\lambda_1\ge\lambda_2\ge\cdots\ge\lambda_{\ell-1}\ge|\lambda_\ell|$ in type $D_\ell$, force every remaining $\lambda_j$ to vanish. The same coordinate argument applies in type $D_3\cong A_3$. Hence $\lambda=\varepsilon_1$, the highest weight of the vector representation, and $\dim_{\mathbb R}V=n+2$.

Finally, $\dim\mathfrak h=(n+2)(n+1)/2=\dim\mathfrak{so}(V,B)$, so the faithful homomorphism $\mathfrak h\to\mathfrak{so}(V,B)$ is an isomorphism. Since the real rank of an orthogonal Lie algebra is the smaller inertia index of its defining form and $\rank_{\mathbb R}\mathfrak h=2$, the smaller inertia index of $B$ is two.
\end{proof}

\subsection{Proof of Theorem~\ref{thm:intro-smallwall}}

\begin{proof}[Proof of Theorem~\ref{thm:classification}]
The rank-one assertion is Proposition~\ref{prop:rankone-small-wall-complete}; it includes the real, complex, and quaternionic types and all low-rank exceptions. The rank-two assertion is Proposition~\ref{prop:ranktwo-smallwall}. The $\Sp_{1,1}\cong\Spin_{4,1}$ vector module already occurs in the orthogonal rank-one row.
\end{proof}

\begin{proof}[Proof of Theorem~\ref{thm:intro-smallwall}]
Suppose that $\limcone_H^{\SO_{p,3}}\subset F_2$. By Proposition~\ref{prop:rankred}, the Lie algebra of $H$ is noncompact simple and has real rank at most two. Moreover, \eqref{eq:small-wall-count} gives $n_V(A)\le2$, so Theorem~\ref{thm:classification} applies. All the invariant forms in its list have smaller inertia index at most two, except for the adjoint module of $\SO_{3,1}^{\circ}$, whose signature is $(3,3)$. None therefore has signature $(p,3)$ with $p>3$, a contradiction.
\end{proof}

\section{Proof of the index-three rigidity theorem}\label{sec:proofmain}

\begin{proof}[Proof of Theorem~\ref{thm:main}]
Let $\Lambda=\Lambda_\rho$. By Proposition~\ref{prop:constant}, the radical rank is constant and equal to $0$, $1$, or $2$.

If $r=2$, Proposition~\ref{prop:maxradpq} gives a contradiction. If $r=1$, Lemma~\ref{lem:spectral} gives the wall condition $\lambda_3=0$ for every infinite-order element. If $r=0$, Lemma~\ref{lem:spectral} gives the stronger condition $\lambda_2=\lambda_3=0$, hence the same wall condition $\lambda_3=0$ a fortiori.

Let $\Delta=\rho(\Gamma)<\PO_{p,3}$, let $\widetilde\Delta<\Ogrp_{p,3}$ be its full inverse image, and put $H=\Zcl(\widetilde\Delta)$. By Corollary~\ref{cor:red}, $H^\circ$ is reductive and $H$ has no degenerate invariant subspace. Suppose first that $H$ is irreducible, and put $\widetilde\Delta_0:=\widetilde\Delta\cap H^\circ$. This subgroup has finite index in $\widetilde\Delta$, has the same boundary and limit set, and is Zariski dense in $H^\circ$; every infinite-order element acts proximally on $\Pj(V)$, with the same continuous attracting-line map. We claim that $H^\circ$ is irreducible. Otherwise, complete reducibility gives an $H^\circ$-invariant direct-sum decomposition. Lemma~\ref{lem:directsum}, applied to $\widetilde\Delta_0$, puts $\Lambda$ in one proper summand. But $\operatorname{span}(\Lambda)$ is $H$-invariant and nonzero, so irreducibility of $H$ gives $\operatorname{span}(\Lambda)=V$, a contradiction.

Lemma~\ref{lem:semisimple-irreducible} shows that $H^\circ$ is semisimple. It is noncompact, since a discrete subgroup of a finite extension of a compact group would be finite. The cone $\limcone_{\widetilde\Delta_0}^{\SO_{p,3}}$ is contained in $F_2$, since passing to a finite-index subgroup does not change the limit cone. Lemma~\ref{lem:wall-inheritance}, applied to $L=H^\circ$, $\Delta=\widetilde\Delta_0$, and $m=2$, gives $\limcone_{H^\circ}^{\SO_{p,3}}\subset F_2$. Theorem~\ref{thm:intro-smallwall} now gives a contradiction. Therefore $H$ is reducible. Since here $q=3$, Proposition~\ref{prop:descent} gives a preserved nondegenerate subspace of signature $(p,b)$ with $b\le2$; for $b=2$ the Glorieux--Monclair theorem \cite[Theorem~1.3]{GM24} shows that the projected representation in $\PO_{p,2}$ is Fuchsian, and for $b=1$ it already preserves a totally geodesic $\HH^p$.

Thus in all cases $\rho(\Gamma)$ preserves a totally geodesic $\HH^p$.

Let $\pi:\Delta\longrightarrow\PO_{p,1}$ be projection along the compact $\Ogrp_2$-factor. Since $\Delta$ is discrete and the second factor is compact, $\pi(\Delta)$ is discrete and $F=\ker\pi$ is finite. As in the proof of Proposition~\ref{prop:descent}, projection leaves the first Jordan gap unchanged, so $\pi\circ\rho$ is $P_1$-Anosov. The rank-one characterization of convex cocompact representations therefore shows that $\pi(\Delta)$ is convex cocompact.

Let $X$ denote the preserved copy of $\HH^p$. The compact complementary factor has all eigenvalues of modulus one, so every infinite-order element of $\Delta$ and its image under $\pi$ have the same attracting eigenline, lying in $\partial_\infty X$. Taking closures of these attracting points shows that $\pi(\Delta)$ has limit set $\Lambda\subset\partial_\infty X$. Both $\Lambda$ and $\partial_\infty X$ are homeomorphic to $\mathbb S^{p-1}$. Invariance of domain makes $\Lambda$ open in $\partial_\infty X$, and compactness makes it closed. Hence $\Lambda=\partial_\infty X$, and convex cocompactness implies that $\pi(\Delta)$ is a uniform lattice.

Finally, Selberg's lemma \cite{Sel60} gives a torsion-free subgroup $\Delta'<\Delta$ of finite index. Since $F$ is finite, $\Delta'\cap F=\{1\}$. Projection therefore identifies $\Delta'$ with a finite-index subgroup of the uniform lattice $\pi(\Delta)$. Thus $\Gamma/\ker\rho\cong\Delta$ is virtually a uniform lattice in $\operatorname{Isom}(\HH^p)$.
\end{proof}

\section{General signatures and proof of \texorpdfstring{Theorem~\ref{thm:oddpq}}{the odd-p theorem}}\label{sec:pqodd}

We return to the general-signature setting of Subsections~\ref{sec:graph} and~\ref{sec:tangent}. Let $p>q\ge2$. By Proposition~\ref{prop:constant}, the tangent radical is a continuous vector bundle of constant rank $0\le r\le q-1$. The maximal rank $r=q-1$ has already been excluded by Proposition~\ref{prop:maxradpq}. We now show that when $p$ is odd, no nonzero radical rank can occur.

\begin{proposition}[Exclusion of nonzero radical rank for odd $p$]\label{prop:oddpqrad}
Under the hypotheses of Theorem~\ref{thm:oddpq}, the tangent radical has rank zero. In particular, $\Lambda_\rho$ is spacelike.
\end{proposition}
\begin{proof}
Put $n=p-1$. Then $n$ is even and $\mathbb S^n$ is simply connected with $\chi(\mathbb S^n)=2$. By Proposition~\ref{prop:constant}, $R=\Rad(Q)$ is a continuous rank-$r$ subbundle of $T\mathbb S^n$. Suppose $0<r\le q-1$. Since $p>q$, we have $0<r<n$. The bundle $R$ and its orthogonal complement $R^\perp$, taken with respect to an auxiliary Riemannian metric, are orientable because $H^1(\mathbb S^n;\Z/2)=0$. Choose their orientations compatibly with the splitting $T\mathbb S^n=R\oplus R^\perp$. The Whitney-sum formula for the Euler class, as stated in \cite[Property~9.6]{MilnorStasheff}, gives
$$
e(T\mathbb S^n)=e(R)\smile e(R^\perp).
$$

But $e(R)\in H^r(\mathbb S^n;\mathbb Z)$ and $e(R^\perp)\in H^{n-r}(\mathbb S^n;\mathbb Z)$, and both groups vanish because $0<r<n$. Hence $e(T\mathbb S^n)=0$, contradicting
$$
\int_{\mathbb S^n} e(T\mathbb S^n)=\chi(\mathbb S^n)=2.
$$

Here we use the equality between the Euler number and the Euler characteristic; see \cite[Corollary~11.12]{MilnorStasheff}. Thus, no nonzero radical rank occurs.
\end{proof}

\begin{proof}[Proof of Theorem~\ref{thm:oddpq}]
By Proposition~\ref{prop:oddpqrad}, the limit set is spacelike. Proposition~\ref{prop:r0entropy} now shows that $\rho(\Gamma)$ preserves a totally geodesic subspace $X$ isometric to $\HH^p$. Equivalently, after conjugacy, $\rho(\Gamma)$ is contained in the image of $\Ogrp_{p,1}\times\Ogrp_{q-1}$ under the natural block-diagonal homomorphism to $\PO_{p,q}$.

Let $\Delta=\rho(\Gamma)$ and let $\pi:\Delta\to\PO_{p,1}$ be the restriction of the projection along the compact $\Ogrp_{q-1}$-factor. As in the final part of the proof of Theorem~\ref{thm:main}, $\ker\pi$ is finite and $\pi(\Delta)$ is a discrete convex cocompact subgroup of $\PO_{p,1}$ with limit set $\Lambda$. Since $\Lambda$ and $\partial_\infty X$ are embedded $(p-1)$-spheres and $\Lambda\subset\partial_\infty X$, invariance of domain gives $\Lambda=\partial_\infty X$. The group $\pi(\Delta)$ is therefore a uniform lattice. Applying Selberg's lemma as in the proof of Theorem~\ref{thm:main} shows that $\Gamma/\ker\rho\cong\Delta$ is virtually a uniform lattice in $\operatorname{Isom}(\HH^p)$.
\end{proof}

\begin{remark}
For even $p$, the Euler-class argument gives no obstruction, since $\chi(\mathbb S^{p-1})=0$. Proposition~\ref{prop:maxradpq} still excludes the maximal radical rank $r=q-1$. The intermediate ranks $1\le r\le q-2$ require additional arguments, such as those used in Section~\ref{sec:proofmain} for $q=3$.
\end{remark}

\section{Further questions}
Conjecture~\ref{conj:pq} remains open when $p$ is even and $q\ge4$. The spacelike case is settled by Proposition~\ref{prop:r0entropy}, the maximal radical rank is excluded by Proposition~\ref{prop:maxradpq}, and a reducible Zariski closure allows restriction to a subspace with smaller negative index by Proposition~\ref{prop:descent}. For an intermediate radical rank $1\le r\le q-2$, Corollary~\ref{cor:wall} places the limit cone in $F_{r+1}$, defined by $\lambda_{r+2}=0$. If the connected Zariski closure is semisimple and acts irreducibly, Lemma~\ref{lem:wall-inheritance} places its ambient limit cone in the same face and bounds its real rank by $r+1$. The remaining question is which groups satisfying these restrictions can occur as Zariski closures of representations with the prescribed limit sphere. Remark~\ref{rem:section-six-general-q} shows that containment in the last wall alone does not exclude irreducible groups in higher signatures.

For representations into $\PO_{p,3}$, a different possible route to Theorem~\ref{thm:main} is analytic. By Seppi--Smith--Toulisse \cite[Corollary~1.1]{SST}, a convex cocompact representation has a unique invariant complete maximal spacelike $p$-submanifold $M_\rho\subset\HH^{p,2}$ with ideal boundary $\Lambda_\rho$. Mazzoli--Viaggi \cite[Proposition~2.7]{MV26} show that $M_\rho$ is totally geodesic if every pseudo-Riemannian Busemann function centered at a point of $\Lambda_\rho$ restricts to a function on $M_\rho$ whose gradient has norm at most $1$. It would therefore suffice to derive this bound directly from the $C^1$ regularity of $\Lambda_\rho$. One may compare this question with the compactness and degeneration arguments for maximal surfaces in \cite[Proposition~2.28 and Sections~4.4--4.7, especially Propositions~4.12 and~4.16]{LT}.

\bibliographystyle{siam}

\bibliography{biblio.bib}

\end{document}